\documentclass[11pt]{article}

\usepackage{amsmath, amsfonts, amssymb,latexsym,wasysym,graphicx}
\input epsf.sty

\newtheorem{Theorem}{Theorem}[section]
\newtheorem{Lemma}{Lemma}[section]
\newtheorem{Proposition}[Lemma]{Proposition}
\newtheorem{Corollary}[Lemma]{Corollary}
\newtheorem{Example}[Lemma]{Example}
\newtheorem{Definition}[Lemma]{Definition}

\newcommand{\BEQ}{\begin{equation}}     % Gleichungen Anfang ..
\newcommand{\BEA}{\begin{eqnarray}}
\newcommand{\BD}{\begin{displaymath}}
\newcommand{\EEQ}{\end{equation}}       % .. und Ende
\newcommand{\EEA}{\end{eqnarray}}
\newcommand{\ED}{\end{displaymath}}
\newcommand{\del}{\delta}
\newcommand{\Del}{\Delta}
\newcommand{\eps}{\varepsilon}          % epsilon

\newcommand{\tgoth}{{\mathfrak{t}}}

\newcommand{\R}{\mathbb{R}}

\newcommand{\Z}{\mathbb{Z}}
\newcommand{\N}{\mathbb{N}}
\newcommand{\I}{{\mathbb{I}}}

\newcommand{\SkI}{{\mathrm{SkI}}}
\newcommand{\supp}{{\mathrm{supp}}}

\def\T{{\mathbb{T}}}

\def\sgn{{\mathrm{sgn}}}

\def\shuffle{\pitchfork}

\newcommand{\eop}{\hfill $\Box$}        % quod erat demonstrandum ... 

\newcommand{\II}{{\rm i}}               % gerades i fuer komplexe Einheit
\newcommand{\half}{{1\over 2}}          % 1/2 als Bruch
\renewcommand{\vec}[1]{\boldsymbol{#1}} % Vektoren fettgedruckt
\catcode`\@=11
\def\numberbysection{\@addtoreset{equation}{section}
        \def\theequation{\thesection.\arabic{equation}}}
\numberbysection

\begin{document}

\vspace*{1.5cm}
\begin{center}
{\Large \bf H\"older-continuous rough paths by Fourier normal ordering
}
\end{center}

\vspace{2mm}
\begin{center}
{\bf  J\'er\'emie Unterberger}
\end{center}

\vspace{2mm}
\begin{quote}

\renewcommand{\baselinestretch}{1.0}
\footnotesize
{We construct in this article an explicit geometric rough path over arbitrary
$d$-dimensional paths with finite $1/\alpha$-variation
 for any $\alpha\in(0,1)$. The method 
 may be coined as 'Fourier normal ordering', since it consists in a regularization obtained
after permuting the order of integration
in iterated integrals so that innermost integrals have highest  Fourier frequencies. In doing so,
there appear non-trivial tree combinatorics, which are best understood by using the structure
of the  Hopf algebra  of decorated rooted trees (in connection with the Chen or multiplicative property) and
of the Hopf shuffle algebra (in connection with the shuffle or geometric property).
 H\"older continuity is proved by using  Besov norms.

The method is well-suited in particular in view of applications to probability theory (see the companion article
\cite{Unt09} for the construction of a rough path over multidimensional fractional Brownian motion
with Hurst index $\alpha<1/4$, or \cite{Unt09ter} for a short survey in that case).
 }
\end{quote}

\vspace{4mm}


\tableofcontents

\newpage

%%%%%%%%%%%%%%%%%%%%%%%%%%%%%%%%%%%%%%%%%%%%%%%%%%%%%%%%%%%%%5
%%%%%%%%%%%%%%%%%%%%%%%%%%%%%%%%%%%%%%%%%%%%%%%%%%%%%%%%%%%%%%%%5

\section{Introduction}

%%%%%%%%%%%%%%%%%%%%%%%%%%%%%%%%%%%%%%%%%%%%%%%%%%%%%%%%%%%%%%5
%%%%%%%%%%%%%%%%%%%%%%%%%%%%%%%%%%%%%%%%%%%%%%%%%%%%%%%%%%%%%%%%

 Assume
$t\mapsto \Gamma_t=(\Gamma_t(1),\ldots,\Gamma_t(d))$, $t\in\R$
 is a smooth $d$-dimensional path, and let $V_1,\ldots,V_d:\R^d\to\R^d$
be smooth vector fields. Then the classical Cauchy-Lipschitz theorem implies that the differential
equation driven by $\Gamma$
\BEQ dy(t)=\sum_{i=1}^d V_i(y(t))d\Gamma_t(i) \label{eq:0:eq-dif} \EEQ
admits a unique solution with initial condition $y(0)=y_0$. The usual way to prove this is
to show by a functional  fixed-point theorem that iterated integrals 
\BEQ y_n\mapsto
y_{n+1}(t):=y_0+\int_0^t \sum_i V_i(y_n(s))d\Gamma_s(i) \EEQ
converge when $n\to \infty$.

Assume now that $\Gamma$ is only 
 $\alpha$-H\"older continuous for some $\alpha\in(0,1)$. Then the Cauchy-Lipschitz theorem
does not hold any more because one first needs to give a meaning to the above integrals, and
in particular to the iterated integrals
$\int_s^t d\Gamma_{t_1}(i_1)\int_s^{t_1} d\Gamma_{t_2}(i_2)
\ldots \int_{s}^{t_{n-1}} d\Gamma_{t_n}(i_n),$ $n\ge 2$, $1\le i_1,\ldots,i_n\le d$.

\medskip

The theory of rough paths, invented by T. Lyons \cite{LQ-bk} and further developed by V. Friz,
N. Victoir \cite{FV} and M. Gubinelli \cite{Gu} implies the possibility to solve eq. (\ref{eq:0:eq-dif})
by a {\em redefinition of the integration} along $\Gamma$, using as an essential ingredient
a {\em rough path} $\bf\Gamma$ over $\Gamma$. By definition, a functional ${\bf\Gamma}=(\vec{\Gamma}^{1},\ldots,\vec{\Gamma}^{N})$, 
 $N=\lfloor 1/\alpha \rfloor$=entire part of $1/\alpha$, is called a {\em rough path over $\Gamma$} if
  $\vec{\Gamma}^{1}_{ts}=(\del\Gamma)_{ts}:=\Gamma_t-\Gamma_s$ are the two-point
increments of $\Gamma$, and
$\vec{\Gamma}^k=(\vec{\Gamma}^k(i_1,\ldots,i_k))_{1\le i_1,\ldots,i_k\le d}$, $k=1,\ldots,N$ 
satisfy the following three properties:

\begin{itemize}
\item[(i)] (\it H\"older continuity)
 each component of $\vec{\Gamma}^k$, $k=1,\ldots,N$
 is 
$k\alpha$-H\"older continuous, that is
to say, $\sup_{s\in\R} \left(\sup_{t\in\R} \frac{|{\bf \Gamma}_{ts}^k(i_1,\ldots,i_k)|}{|t-s|^{k\alpha}} \right)<\infty.$

\item[(ii)] ({\it multiplicative/Chen property}) letting $\del{\bf\Gamma}^k_{tus}:=
{\bf\Gamma}_{ts}^k-{\bf\Gamma}^k_{tu}-{\bf \Gamma}^k_{us}$, one requires
\BEQ
 \del\vec{\Gamma}^k_{tus}(i_1,\ldots,i_k) = \sum_{k_1+k_2=k} \vec{\Gamma}_{tu}^{k_1}(i_1,\ldots,i_{k_1}) \vec{\Gamma}_{us}^{k_2}(i_{k_1+1},\ldots,i_k); \label{eq:0:x}
\EEQ

\item[(iii)] ({\it geometric/shuffle property}) 
\BEQ  {\bf \Gamma}^{n_1}_{ts}(i_1,\ldots,i_{n_1}) {\bf \Gamma}^{n_2}_{ts}(j_1,\ldots,j_{n_2}) 
 = 
\sum_{\vec{k}\in {\mathrm{Sh}}(\vec{i},\vec{j})} {\bf \Gamma}^{n_1+n_2}(k_1,\ldots,k_{n_1+n_2})  \label{eq:0:geo} \EEQ
where ${\mathrm{Sh}}(\vec{i},\vec{j})$ is the set of {\em shuffles} of $\vec{i}=(i_1,\ldots,i_{n_1})$
and $\vec{j}=(j_1,\ldots,j_{n_2})$, that is to say, of  permutations of $i_1,\ldots,i_{n_1},j_1,\ldots,j_{n_2}$
which do not change the orderings of $(i_1,\ldots,i_{n_1})$ and $(j_1,\ldots,j_{n_2})$.
\end{itemize}

\medskip

There is a canonical choice for $\bf\Gamma$, called {\em canonical lift of $\Gamma$},
 when $\Gamma$ is a smooth path, namely, the iterated integrals of $\Gamma$ of arbitrary order.
If one sets 
\BEQ {\bf\Gamma}^{cano,n}(i_1,\ldots,i_n):=\int_s^t d\Gamma_{t_1}(i_1)\int_s^{t_1} d\Gamma_{t_2}(i_2)
\ldots \int_{s}^{t_{n-1}} d\Gamma_{t_n}(i_n), \EEQ
then ${\bf\Gamma}^{cano}=({\bf\Gamma}^{cano})_{n=1,2,\ldots}$ satisfies properties (i), (ii), (iii)
with $\alpha=1$. Axiom (ii) receives a natural geometric interpretation in this case since ${\bf\Gamma}^{cano}$
measures the areas, volumes and so forth generated by $\Gamma^1,\ldots,\Gamma^d$, see \cite{FV}, while
axiom (iii) may be deduced from Fubini's theorem. A further justification of axioms (i),(ii),(iii)
comes from the fact that any rough path is a {\em limit in some sense of the iterated integrals of a 
sequence of smooth paths}, so ${\bf\Gamma}$ plays the r\^ole of a {\em substitute} of iterated integrals
for $\Gamma$.

\bigskip

The problem we address here is the {\em existence and construction} of rough paths. 
It is particularly relevant when $\Gamma$ is a random path; it  allows for  the pathwise
construction of stochastic integrals or of solutions of stochastic differential equations driven by $\Gamma$.
Rough paths are then usually constructed by choosing some appropriate smooth approximation
$\Gamma^{\eta}$, $\eta\overset{>}{\to} 0$ of $\Gamma$ and proving that the canonical lift
of $\Gamma^{\eta}$
converges in $L^2(\Omega)$  for appropriate H\"older norms to a rough path $\bf \Gamma$ lying above $\Gamma$
(see \cite{CQ02,Unt08} in the case of fractional Brownian motion with Hurst index $\alpha>1/4$, 
and \cite{BHL,HL} for a class of random paths on fractals, or references in \cite{LyoVic07}).

A general construction of a rough path for deterministic paths has been given -- in the original formulation due to
T. Lyons -- in an article by T. Lyons and N. Victoir \cite{LyoVic07}. The idea 
\cite{FV}
is to see a rough path over $\Gamma$ as a H\"older section of the trivial
 $G$-principal bundle over $\R$, where $G$ is a free rank-$N$ nilpotent group (or Carnot group), while
the underlying path $\Gamma$ is a section of the corresponding quotient 
$G/K$-bundle for some normal subgroup $K$ of $G$; so
one is reduced to the problem of finding  H\"older-continuous sections $g_t K\to g_t$. Obviously, 
there is
no canonical way to do this in general. This
abstract, group-theoretic construction  -- which uses the axiom of choice --
 is unfortunately  not particularly
appropriate for concrete problems, such as the behaviour of solutions of stochastic differential
equations for instance.

\bigskip

We propose here a new, {\em explicit method} to construct a rough path $\bf\Gamma$ over an arbitrary
$\alpha$-H\"older path $\Gamma$ which rests on an {\em algorithm} that we call {\em Fourier normal
ordering}. Let us explain the main points of this algorithm. The first point is the use of
{\em Fourier transform}, $\cal F$; H\"older estimates are obtained by means of Besov norms involving
compactly supported Fourier multipliers, see Appendix. Assume for simplicity that $\Gamma$ is
complactly supported; this assumption is essentially void since one may multiply any $\alpha$-H\"older
path by a smooth, compactly supported function equal to $1$ over an arbitrary large compact interval,
and then restrict the construction to this interval.
What makes the Fourier transform interesting for our problem is that $({\cal F}\Gamma')(\xi)=\II\xi
({\cal F}\Gamma)(\xi)$ is a well-defined function; thus, the meaningless iterated integral
$\int_s^t d\Gamma_{t_1}(i_1)\int_s^{t_1} d\Gamma_{t_2}(i_2)
\ldots \int_{s}^{t_{n-1}} d\Gamma_{t_n}(i_n)$ is rewritten after Fourier transformation
as some integral $\int_{-\infty}^{\infty}
\ldots\int_{-\infty}^{+\infty} f(\xi_1,\ldots,\xi_n)d\xi_1\ldots d\xi_n$, where $f$ is regular but not
integrable at infinity along certain directions.

The second, main point is the {\em splitting of the Fourier domain of integration} $\R^n$ into
$\cup_{\sigma\in\Sigma_n} \R^n_{\sigma}$, $\Sigma_n$=set of permutations of $\{1,\ldots,n\}$, where
$\R^n_{\sigma}:=\{|\xi_{\sigma(1)}|\le\ldots\le |\xi_{\sigma(n)}|\}$, see section 2  for a more
accurate definition involving the Besov dyadic decomposition. Away from the singular directions,
the resulting integrals are naturally shown to have a polynomially decreasing behaviour at infinity
implying the correct H\"older behaviour; simple examples may be read from \cite{Unt09ter}. However -- as
computations in section 4 clearly show, see also \cite{Unt09ter} for an elementary example --
these bounds are naturally obtained only {\em after permuting the order of integration} by means
 of Fubini's theorem, 
so that the Fourier coordinates $|\xi_1|,\ldots,|\xi_n|$ {\em appear in increasing order}. There appear
in the process integrals over domains which differ from the simplex $\{t\ge t_1\ge\ldots\ge t_n\ge s\}$,
which are particular instances of {\em tree integrals}, and that we call {\em tree skeleton integrals}.

The next step is to {\em regularize} the tree skeleton integrals so that Fourier integrals converge
at infinity, {\em without} losing the Chen and shuffle properties (ii) and (iii). At this point it turns
out to be both natural and necessary to re-interpret the above scheme in terms of tree  Hopf algebra
combinatorics.
 The interest for the study of Hopf algebras of
trees or graphs surged out of a series of papers by A. Connes and D. Kreimer
\cite{ConKre98,ConKre00,ConKre01} concerning the mathematical structures hidden
behind the Bogolioubov-Hepp-Parasiuk-Zimmermann (BPHZ)
 procedure for renormalizing Feynmann diagrams in quantum field theory \cite{Hepp},
and is still very much alive, see for instance \cite{Kre99,BroFra03,BroFraKra06,Foi02,CEM08,Mur2,Wal00,ChaLiv07},
with applications ranging from numerical methods to quantum chromodynamics or 
multi-zeta functions or operads. It appears that the shuffle property may be stated by saying that {\em regularized
skeleton integrals define characters} of yet another Hopf algebra called {\em shuffle algebra},
while the Chen property follows from the very definition of the regularized iterated
integrals as a {\em convolution of regularized skeleton integrals}. 

We show that the tree skeleton integrals may be regularized by integrating over appropriate
subdomains of $\R^n_{\sigma}$ avoiding the singular directions. The proof of properties (ii), (iii)
uses Hopf combinatorics and does not depend on the choice of the above subdomains, while the
proof of the H\"older estimates (i) uses both tree combinatorics and some elementary analysis relying
on the shape of the subdomains.

\medskip

It seems natural to look for a less arbitrary regularization scheme for the skeleton integrals. The
idea of cancelling singularities by building iteratively counterterms, originated from the BPHZ procedure,
should also apply here. We plan to give such a construction (such as dimensional regularization for
instance) in the near future.

\bigskip

Let us state our main result. Throughout the paper $\alpha\in(0,1)$ is some fixed constant
and $N=\lfloor 1/\alpha\rfloor$.

\medskip

{\bf Main theorem.}

{\it Assume $1/\alpha\not\in\N$. Let $\Gamma=(\Gamma(1),\ldots,\Gamma(d)):\R\to\R^d$ be a compactly supported $\alpha$-H\"older path. Then the functional $({\cal R}{\bf\Gamma}^1,\ldots,
{\cal R}{\bf\Gamma}^N)$ defined in section 2  is an $\alpha$-H\"older  geometric rough path lying over $\Gamma$ in the sense
of properties (i),(ii),(iii) of the Introduction.
}

\medskip

In a companion paper \cite{Unt09},
 we construct by  the same algorithm an explicit rough path over a $d$-dimensional fractional Brownian
motion $B^{\alpha} =(B^{\alpha}(1),\ldots,B^{\alpha}(d))$ with arbitrary Hurst index $\alpha\in(0,1)$ --
recall simply that the paths of $B^{\alpha}$ are a.s. $\kappa$-H\"older for every $\kappa<\alpha$. The
problem was up to  now open for $\alpha\le 1/4$ despite many attempts \cite{CQ02,Unt08,Unt08b,DNN}.
 Fourier normal ordering  turns out to be  very efficient in combination with Gaussian tools, and
provides explicit bounds for the moments of the  rough path, seen as a path-valued random
variable.

\medskip

The above theorem extends to paths $\Gamma$ with finite $1/\alpha$-variation. Namely (see
\cite{LyoVic07}, \cite{Lej03}  or also \cite{FV}),
 a simple change of variable $\Gamma\to \Gamma^{\phi}:=
\Gamma\circ\phi^{-1}$ turns $\Gamma$ into an $\alpha$-H\"older path, with $\phi$ defined for instance
as $\phi(t):=\sup_{n\ge 1}\sup_{0=t_0\le \ldots\le t_n=t} \sum_{j=0}^{n-1}
||\Gamma(t_{j+1})-\Gamma(t_j)||^{1/\alpha}$. The construction of the above Theorem, applied
to $\Gamma^{\phi}$,  yields a family
of paths with H\"older regularities $\alpha,2\alpha,\ldots,N\alpha$ which may alternatively be seen
as a $G^N$-valued $\alpha$-H\"older path ${\bf\Gamma}^{\phi}$,
 where $G^N$ is the Carnot free nilpotent group of order $N$
equipped with any subadditive homogeneous norm. Then (as proved in \cite{LyoVic07}, Lemma 8)
${\bf\Gamma}:={\bf\Gamma}^{\phi}\circ\phi$ has finite $1/\alpha$-variation, which is equivalent to saying that
${\bf\Gamma}^n$ has finite $1/n\alpha$-variation for $n=1,\ldots,N$, and lies above $\Gamma$.

\medskip

{\bf Corollary.}

{\em Let $\alpha\in(0,1)$ and  $\alpha'<\alpha$.
 Then every $\alpha$-H\"older path $\Gamma$ may be lifted to a {\em strong} $\alpha'$-H\"older
 geometric rough path, namely,  there exists
a sequence of canonical lifts ${\bf \Gamma}^{(n)}$ of smooth paths $\Gamma^{(n)}$
 converging to ${\cal R}{\bf\Gamma}$ for the sequence of $\alpha'$-H\"older norms.
}

\medskip

The set of {\it strong} $\alpha$-H\"older geometric rough paths is strictly included in the set of general
$\alpha$-H\"older  geometric rough paths. On the other hand, as we already alluded to  above,
 a weak $\alpha$-H\"older geometric rough path
may be seen as a strong $\alpha'$-H\"older geometric rough path if $\alpha'<\alpha$. This accounts for the
loss of regularity in the Corollary (see \cite{FV} for a precise discussion). The proviso $1/\alpha\not\in\N$
in the statement of the main theorem is a priori needed because otherwise ${\cal R}{\bf \Gamma}^{N}$ may not be
treated in the same way as the lower-order iterated integrals (although we do not know if it is actually
necessary). However, if $1/\alpha\in\N$, all one has to do is replace $\alpha$ by a slightly smaller parameter
$\alpha'$, so that the Corollary holds even in this case.

Note that the present paper gives unfortunately no explicit way of approximating ${\cal R}{\bf \Gamma}$ by
{\em canonical lifts} of smooth  paths, i.e. of seeing it concretely as a {\it strong} geometric rough path. The question is
currently under investigation in the particular case of fractional Brownian motion by using constructive
field theory methods. Interestingly enough, the idea of controlling singularities by separating the
 Fourier scales according to a dyadic decomposition is at the core of constructive field theory 
\cite{Riv}.

\bigskip

Here is an outline of the article. A thorough presentation of iterated integrals, together with the skeleton
integral variant, the implementation of Fourier normal ordering, and the extension to tree integrals,
is given in section 1, where $\Gamma$ is assumed to be {\em smooth}. The {\em regularization algorithm}
is presented in section 2; the regularized rough path ${\cal R}\Gamma$ is defined there for an arbitrary
$\alpha$-H\"older path $\Gamma$. The proof of the Chen and shuffle properties is given in section 3,
where one may also find two abstract but more compact reformulations of the regularization algorithm, see
Lemma \ref{lem:3.5} and Definition \ref{def:alternative}. H\"older estimates are to be found in section 4. Finally, we gathered in an Appendix some technical
facts about Besov spaces required for the construction.

\bigskip

{\bf Notations.} We shall denote by ${\cal F}$ the {\em Fourier transform},
\BEQ {\cal F}:L^2(\R^l)\to L^2(\R^l), f\mapsto {\cal F}(f)(\xi)=\frac{1}{(2\pi)^{l/2}} \int_{\R^l} f(x)e^{-\II \langle
x,\xi\rangle} \ dx. \EEQ
 Throughout the article, $\Gamma:\R\to\R^d$ is some compactly
supported  $\alpha$-H\"older path; sometimes, it is assumed to be smooth.  The {\em permutation group}
 of $\{1,\ldots,n\}$ is denoted by $\Sigma_n$.
 Also, if $a,b:X\to\R_+$ are functions on some set $X$ such that $a(x)\le Cb(x)$ for every
$x\in X$, we shall write $a\lesssim b$. Admissible cuts of a tree $\T$, see subsection 1.3, are usually denoted by $v$
or $w$, and we write $(Roo_{\vec{v}}(\T),Lea_{\vec{v}}(\T))$ ({\em root} part and {\em leaves}) instead of the
traditional notation $(R^c \T,P^c \T)$ due to Connes and Kreimer.

%%%%%%%%%%%%%%%%%%%%%%%%%%%%%%%%%%%%%%%%%%%%%%%%%%
%%%%%%%%%%%%%%%%%%%%%%%%%%%%%%%%%%%%%%%%%%%%%%%%%%%%%%%%

\section{Iterated integrals : smooth case}

%%%%%%%%%%%%%%%%%%%%%%%%%%%%%%%%%%%%%%%%%%%%%%%%%%%%%%%%%
%%%%%%%%%%%%%%%%%%%%%%%%%%%%%%%%%%%%%%%%%%%%%%%%%

Let $t\mapsto \Gamma_t=(\Gamma_t(1),\ldots,\Gamma_t(d))$ be a $d$-dimensional, compactly
supported, {\em smooth} 
path. The purpose of this section is to give proper notations for iterated integrals of  $\Gamma$
and to introduce some tools which will pave the way for  the regularization algorithm.  Subsection 1.1
on tree iterated integrals  is standard, see for instance \cite{ConKre98}. We
introduce {\em permutation graphs} and {\em Fourier normal ordering} for smooth paths in
subsection 1.2. The {\em tree Chen property} -- a generalization of the usual Chen property to tree
iterated integrals -- is recalled in subsection 1.3, in connection with the underlying Hopf algebraic
structure.  Finally, 
a variant of iterated integrals called {\em skeleton integrals} is introduced in subsection 1.4,
together with a variant of the tree Chen property that we call {\em tree skeleton decomposition}.

%%%%%%%%%%%%%%%%%%%%%%%%%%%%%%%%%%%%%%%%%%

\subsection{From iterated integrals to trees}

%%%%%%%%%%%%%%%%%%%%%%%%%%%%%%%%%%%%%%%%%%%%%

It was noted already long time ago \cite{But72} that iterated integrals could be encoded by trees, see
also \cite{Kre99}. This
remark has been exploited in connection with the construction of the rough path
solution of partial, stochastic differential equations  in \cite{Gu2}.
The correspondence between trees and iterated integrals goes simply as follows.

\begin{Definition}

A {\em decorated rooted tree} (to be drawn growing {\em up}) is a finite tree with a distinguished
vertex called {\em root} and edges oriented {\em downwards}, i.e. directed towards the root, such
that every vertex wears a positive integer label called {\em decoration}.

If $\T$ is a decorated rooted tree, we let $V(\T)$ be the set of its vertices (including the root),
and $\ell:V(\T)\to\N$ be its decoration.

\end{Definition}

\begin{Definition} [tree partial ordering]  \label{def:2:connect}

Let $\T$ be a decorated rooted tree.

\begin{itemize}
\item Letting $v,w\in V(\T)$, we say that $v$ {\em connects directly to} $w$, and write
$v\to w$ or equivalently $w=v^-$, if $(v,w)$ is an edge oriented downwards from $v$
to $w$. Note that $v^-$ exists and is unique except if $v$ is the root.

\item If $v_m\to v_{m-1}\to\ldots\to v_1$, then we shall write $v_m\twoheadrightarrow v_1$, and say that
$v_m$ {\em connects to} $v_1$. By definition, all vertices (except the root) connect to the root.

\item Let $(v_1,\ldots,v_{|V(\T)|})$ be an ordering of $V(\T)$. Assume that
$\left( v_i\twoheadrightarrow v_j\right)\Rightarrow\left(i>j\right)$; in particular, $v_1$ is the root.
Then we shall say that the ordering is {\em compatible} with the {\em tree partial ordering} defined
by $\twoheadrightarrow$.
\end{itemize}

\end{Definition}

\begin{figure}[h]
  \centering
   \includegraphics[scale=0.35]{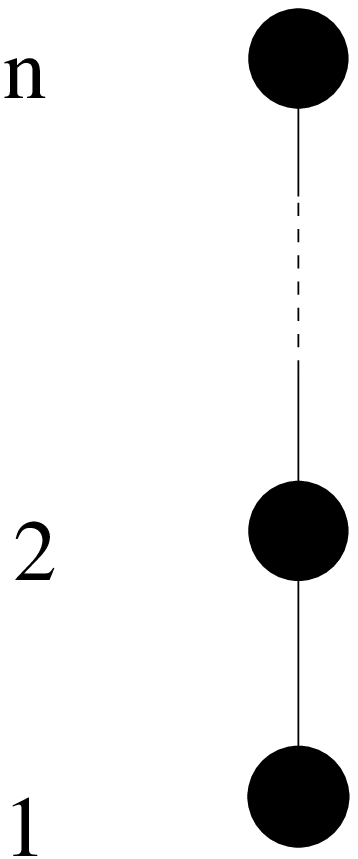}
   \caption{\small{Trunk tree with set of vertices $\{n\to n-1\to \ldots\to 1\}$.}}
  \label{Fig1}
\end{figure}

\begin{Definition}[tree integrals]  \label{def:2:it-int}

\begin{itemize}
\item[(i)]

Let  $\Gamma=(\Gamma(1),\ldots,\Gamma(d))$ be a $d$-dimensional, compactly supported,
 smooth path, and $\T$ a rooted
tree decorated by  $\ell:V(\T)\to\{1,\ldots,d\}$. Then $I_{\T}(\Gamma):\R^2\to\R$ is the iterated integral
defined as
\BEQ [I_{\T}(\Gamma)]_{ts}:=\int_s^t d\Gamma_{x_1}(\ell(v_1))\int_s^{x_{v_2^-}} d\Gamma_{x_2}(\ell(v_2))
\ldots \int_s^{x_{v^-_{|V(\T)|}}} d\Gamma_{x_{v_{|V(\T)|}}}(\ell(v_{|V(\T)|}))  \label{eq:def:ITG}
\EEQ
where $(v_1,\ldots,v_{|V(\T)|})$ is any ordering of $V(\T)$ compatible with the tree partial ordering.

In particular, if $\T$ is a trunk tree with $n$ vertices (see Fig. \ref{Fig1}) -- so that
the tree ordering is total -- we shall write
\BEQ I_{\T}(\Gamma)=I_n^{\ell}(\Gamma),\EEQ
where
\BEQ [I_n^{\ell}(\Gamma)]_{ts}:=\int_s^t d\Gamma_{x_1}(\ell(1)) \int_s^{x_1} d\Gamma_{x_2}(\ell(2))
\ldots \int_s^{x_{n-1}} d\Gamma_{x_{n}}(\ell(n)).
\EEQ

\item[(ii)] (multilinear extension) Assume $\mu$ is a compactly supported, signed Borel measure on $\R^{V(\T)}:=
\{(x_v)_{v\in V(\T)}, x_v\in\R\}$. Then
\BEQ [I_{\T}(\mu)]_{ts}:=\int^t_s \int^{x_{v_2^-}}_s\ldots\int^{x_{v^-_{V(\T)}}}_s \mu(dx_{v_1},\ldots,dx_{v_{V(\T)}}).
\EEQ

\end{itemize}

\end{Definition}

Clearly, the definition of $[I_{\T}(\Gamma)]_{ts}$ given in eq. (\ref{eq:def:ITG}) does not depend on 
the choice of the ordering $(v_1,\ldots,v_{|V(\T)|})$. For instance, consider $\T=\T_1^{\sigma}$ 
to be the first tree in Fig. \ref{Fig4bis}. Then
\BEA  [I_{\T}(\Gamma)]_{ts} &=& \int_s^t d\Gamma_{x_1}(1) \left( \int_s^{x_1} d\Gamma_{x_2}(2)
\int_s^{x_1} d\Gamma_{x_3}(3)\right) \nonumber\\
&=& \int_s^t d\Gamma_{x_1}(1) \left( \int_s^{x_1} d\Gamma_{x_2}(3)
\int_s^{x_1} d\Gamma_{x_3}(2)\right).\EEA

Note that the decoration  of $\T$ is required only for (i). In case of ambiguity, we shall also
use the decoration-independent notation $I_{\T}\left(\otimes_{v\in V(\T)} \Gamma(\ell(v))\right)$ instead
of $I_{\T}(\Gamma)$.

\smallskip

The above correspondence extends by multilinearity to the algebra of decorated rooted trees defined
by Connes and Kreimer \cite{ConKre98}, whose definition  we now recall.

\begin{Definition}[algebra of decorated rooted trees]

\begin{itemize}
\item[(i)] Let $\cal T$ be the set of decorated rooted trees.
\item[(ii)] Let $\bf H$ be the free commutative algebra over $\R$ generated by $\cal T$,
with unit element denoted by $e$.
If $\T_1,\T_2,\ldots \T_l$ are decorated rooted trees, then the product  $\T_1\ldots\T_l$ is the
{\em forest} with connected components $\T_1,\ldots,\T_l$. 
\item[(ii)] Let $\T'=\sum_{l=1}^L m_l \T_l\in{\bf H}$, where $m_l\in\Z$ and each $\T_l=\T_{l,1}\ldots
\T_{l,j_l}$ is a forest whose decorations have values  in the set $\{1,\ldots,d\}$.  Then
\BEQ [I_{\T'}(\Gamma)]_{ts}:=\sum_{l=1}^L m_l [I_{\T_{l,1}}(\Gamma)]_{ts} \ldots [I_{\T_{l,j_l}}(\Gamma)]_{ts}.
 \EEQ
\end{itemize}

\end{Definition}

%%%%%%%%%%%%%%%%%%%%%%%%%%%%%%%

\subsection{Permutation graphs and Fourier normal ordering for smooth paths}

%%%%%%%%%%%%%%%%%%%%%%%%%%%%%%%%%%%%%%%%

As explained briefly in the Introduction, and as we shall see in the next sections, an
essential step in our regularization algorithm is to rewrite iterated integrals by permuting
the order of integration. We shall prove the following lemma in this subsection:

\begin{Lemma}[permutation graphs] \label{lem:2:sigma}

To every trunk tree $\T_n$ with $n$ vertices and decoration $\ell$, and every
 permutation $\sigma\in\Sigma_n$,  is associated in a canonical way
 an element $\T^{\sigma}$ of $\bf H$ called {\em permutation graph},
such that:
\begin{itemize}
\item[(i)] \BEQ I_n^{\ell}(\Gamma)=I_{\T^{\sigma}}(\Gamma); \EEQ
\item[(ii)] 
\BEQ \T^{\sigma}=\sum_{j=1}^{J_{\sigma}} g(\sigma,j) \T_j^{\sigma}\in {\bf H} , \label{eq:gsigmaj} \EEQ
where  $g(\sigma,j)=\pm 1$ and each $\T_j^{\sigma}$, $j=1,\ldots, J_{\sigma}$ is a forest   provided
by construction  with a total ordering compatible with its tree structure, image of the
ordering $\{v_1<\ldots<v_n\}$ of the trunk tree $\T_n$ by the permutation $\sigma$. The
decoration of $\T^{\sigma}$ is $\ell\circ\sigma$. 
\end{itemize}
\end{Lemma}

{\bf Proof.}  Let $\sigma\in\Sigma_n$. Applying Fubini's theorem yields
\BEA [I_n^{\ell}(\Gamma)]_{ts}&=& \int_s^t d\Gamma_{x_1}(\ell(1))\int_s^{x_1} d\Gamma_{x_2}(\ell(2))
\ldots \int_s^{x_{n-1}} d\Gamma_{x_{n}}(\ell(n)) \nonumber\\
&=& \int_{s_1}^{t_1} d\Gamma_{x_{\sigma(1)}}(\ell(\sigma(1)))\int_{s_2}^{t_2} d\Gamma_{x_{\sigma(2)}}
(\ell(\sigma(2))) \ldots  \int_{s_{n}}^{t_{n}} d\Gamma_{x_{\sigma(n)}}(\ell(\sigma(n))), 
\nonumber\\ \label{eq:2:2.8}
\EEA
with $s_1=s$, $t_1=t$, and for some suitable choice of $s_j\in\{s\}\cup\{x_{\sigma(i)}, i<j\}$,\ \ 
 $t_j\in\{t\}\cup\{x_{\sigma(i)},i<j\}$
$(j\ge 2)$. Now decompose $\int_{s_j}^{t_j} d\Gamma_{x_{\sigma(j)}}(\ell(\sigma(j)))$ into 
$$\left( \int_s^{t_j}-\int_s^{s_j}\right) d\Gamma_{x_{\sigma(j)}}(\ell(\sigma(j)))$$ if
$s_j\not=s,t_j\not=t$, and $\int_{s_j}^{t} d\Gamma_{x_{\sigma(j)}}(\ell(\sigma(j)))$ into 
$$\left( \int_s^{t}-\int_s^{s_j}\right) d\Gamma_{x_{\sigma(j)}}(\ell(\sigma(j)))$$ if
$s_j\not=s$. Then $I_n^{\ell}(\Gamma)$ has been rewritten as a sum of terms of the form
\BEQ\pm \int_s^{\tau_1} d\Gamma_{x_1}(\ell(\sigma(1)))\int_s^{\tau_2}d\Gamma_{x_2}(\ell(\sigma(2)))
\ldots\int_s^{\tau_{n}} d\Gamma_{x_{n}}(\ell(\sigma(n))), \label{eq:2:2.9} \EEQ
where $\tau_1=t$ and $\tau_j\in\{t\}\cup\{x_i, i<j\}$, $j=2,\ldots,n$. Note the renaming of 
variables and vertices from eq. (\ref{eq:2:2.8}) to eq. (\ref{eq:2:2.9}).  Encoding each of these
expressions by the forest  $\T$ with set of vertices  $V(\T)=\{1,\ldots,n\}$, label function
$\ell\circ\sigma$,  roots $\{j=1,\ldots,n\ |\
\tau_j=t\}$, and oriented edges $\{(j,j^-)\ |\ j=2,\ldots,n, \tau_j=x_{j^-}\}$, yields
\BEQ I_n^{\ell}(\Gamma)=I_{\T^{\sigma}}(\Gamma) \EEQ
for some $\T^{\sigma}\in{\bf H}$ as in eq. (\ref{eq:gsigmaj}). \hfill \eop

\medskip

\begin{figure}[h]
  \centering
   \includegraphics[scale=0.35]{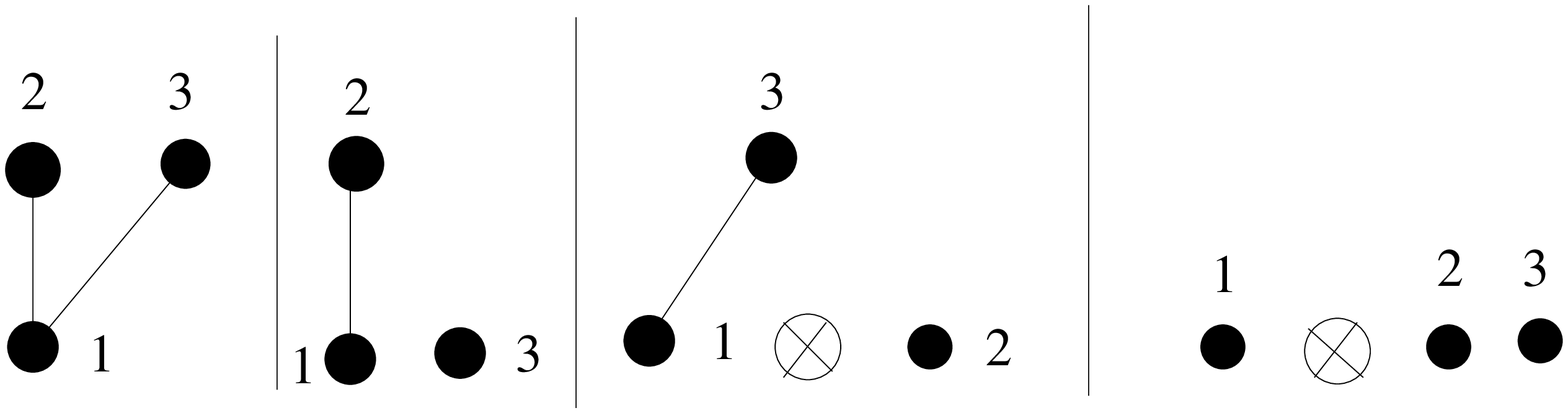}
   \caption{\small{Example \ref{ex:2:1}. From left to right: $\T^{\sigma}_1; \T^{\sigma}_2$;
 $Roo_{\{2\}}\T^{\sigma}_1 \otimes Lea_{\{2\}}\T^{\sigma}_1$; $Roo_{\{2,3\}}\T^{\sigma}_1
\otimes Lea_{\{2,3\}}\T^{\sigma}_1$ }}
 \label{Fig4bis}
\end{figure}

\begin{Example} \label{ex:2:1}

 Let $\sigma=\left(\begin{array}{ccc} 1 & 2 & 3\\ 2 & 3 & 1 \end{array}\right)$.
Then 
\BEA &&  \int_s^t d\Gamma_{x_1}(\ell(1)) \int_s^{t_2} d\Gamma_{x_2}(\ell(2)) \int_s^{t_3}
d\Gamma_{x_3}(\ell(3))= \nonumber\\
&& \qquad -\int_s^t d\Gamma_{x_2}(\ell(2)) \int_s^{x_2} d\Gamma_{x_3}(\ell(3))
\int_s^{x_2} d\Gamma_{x_1}(\ell(1)) \nonumber\\
&& \qquad \qquad  + \int_s^t d\Gamma_{x_2}(\ell(2)) \int_s^{x_2}
d\Gamma_{x_3}(\ell(3)) \ .\ \int_s^t d\Gamma_{x_1}(\ell(1))\EEA
Hence $\T^{\sigma}=-\T_1^{\sigma}+\T_2^{\sigma}$ is the sum of a tree and of a forest with two
components. See Fig. \ref{Fig4bis}, where variables and vertices have been renamed according to the permutation $\sigma$. 

\end{Example}

%%%%%%%%%%%%%%%%%%%%%%%%%%%%%%%

\subsection{Tree Chen property and coproduct structure}

%%%%%%%%%%%%%%%%%%%%%%%

The Chen property (ii), see Introduction, may be generalized to tree iterated integrals by using
the coproduct structure of $\bf H$, as explained in \cite{ConKre98}. It is an essential feature of 
our algorithm since it implies the possibility to reconstruct a rough path $\bf\Gamma$ from the
quantities $t\mapsto {\bf \Gamma}^n_{ts_0}$ with {\em fixed} $s_0$.  This idea will be pursued further in the
next subsection, where  we shall introduce a variant of these iterated integrals with fixed $s_0$
called {\em skeleton integrals}.

\begin{Definition}[admissible cuts] (see \cite{ConKre98}, section 2)

{\it
\begin{enumerate}
\item
 Let $\T$ be a tree, with set of vertices $V(\T)$ and root
denoted by  $0$.
If $\vec{v}=(v_1,\ldots,v_J)$, $J\ge 1$  is any totally disconnected subset of $V(\T)\setminus\{0\}$,
 i.e. $v_i\not\twoheadrightarrow
v_j$ for all $i,j=1,\ldots,J$, then we shall say that $\vec{v}$ is an {\em admissible cut} of $\T$, and
 write $\vec{v}\models V(\T)$. We let $Lea_{\vec{v}}\T$ (read: {\em leaves} of $\T$) 
be the sub-forest (or sub-tree if $J=1$) obtained by keeping
only the vertices above $\vec{v}$, i.e. $V(Lea_{\vec{v}}\T)=\vec{v}\cup\{w\in V(\T):\ \exists j=1,\ldots,J,
w\twoheadrightarrow v_j\}$, and $Roo_{\vec{v}}\T$ (read: {\em root part} of $\T$)
 be the sub-tree obtained by keeping all other vertices.

\item Let $\T=\T_1\ldots\T_l$ be a forest, together with its decomposition into trees. Then an {\em
admissible
cut} of $\T$ is a disjoint union $\vec{v}_1\cup\ldots\cup\vec{v}_l$, $\vec{v}_i\subset\T_i$, where $\vec{v}_i$
is either $\emptyset$, $\{0_i\}$ (root of $\T_i$) or an admissible cut of $\T_i$; by convention,
 the two trivial cuts $\emptyset\cup\ldots\cup\emptyset $ and $\{0_1\}\cup\ldots\cup\{0_l\}$ are
excluded.
 By
definition, we let $Roo_{\vec{v}}\T=Roo_{\vec{v}_1}\T_1 \ldots Roo_{\vec{v}_l}\T_l$, $Lea_{\vec{v}}\T=
Lea_{\vec{v}_1}\T_1\ldots Lea_{\vec{v}_l}\T_l$ (if $\vec{v}_i=\emptyset$, resp. $\{0_i\}$, then
$(Roo_{\vec{v}_i}\T_i,Lea_{\vec{v}_i}\T_i):=(\T_i,\emptyset)$, resp. $(\emptyset,\T_i)$). 

\end{enumerate}
}
\label{def:6:admissible-cut}
\end{Definition}

\begin{figure}[h]
  \centering
   \includegraphics[scale=0.35]{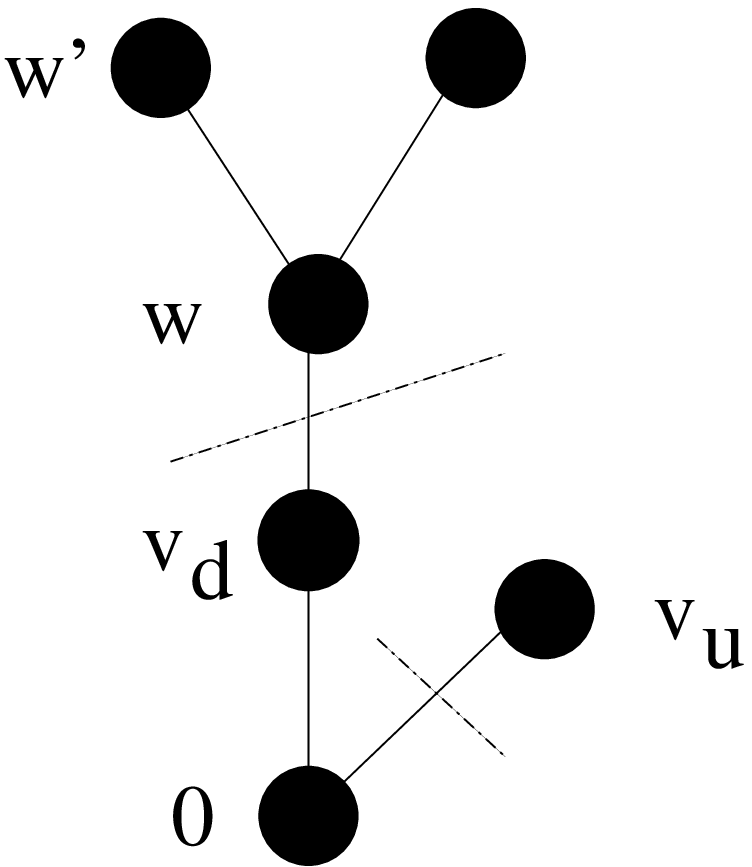}
   \caption{\small{Admissible cut.}}
  \label{Fig2}
\end{figure}

\begin{figure}[h]
  \centering
   \includegraphics[scale=0.35]{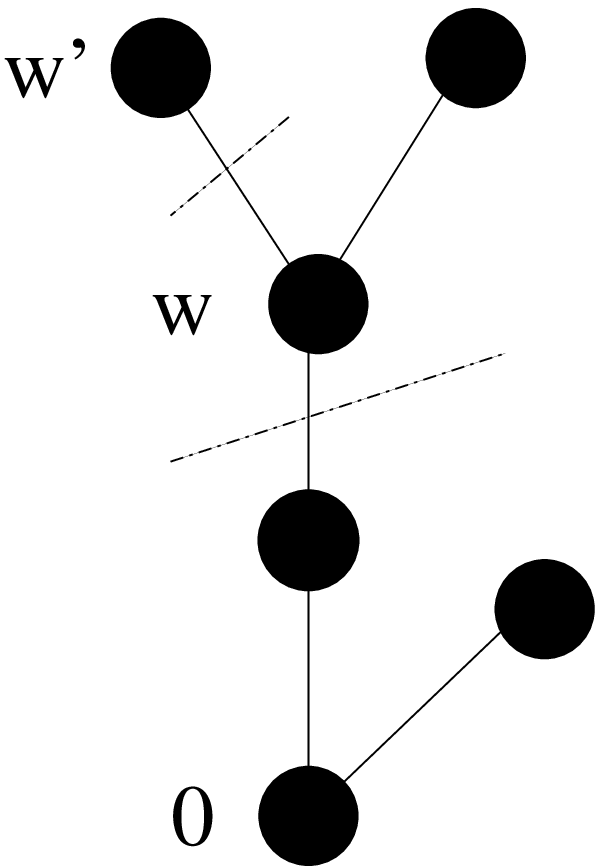}
   \caption{\small{Non-admissible cut.}}
  \label{Fig3}
\end{figure}

See Fig. \ref{Fig2}, \ref{Fig3} and \ref{Fig4bis}. Defining the co-product operation 
\BEA  \Del: {\bf H}\to {\bf H}\otimes
{\bf H}, \quad  \T\mapsto e \otimes \T+\T\otimes e+\sum_{\vec{v}\models V(\T)} Roo_{\vec{v}}\T\otimes
Lea_{\vec{v}}\T \nonumber\\ \EEA
 where $e$ stands for the unit element,
yields a coalgebra structure on ${\bf H}$. One may also define an antipode $S$,
 which makes $\bf H$ a Hopf algebra (see section 3 for more details).

We may now state the {\em tree Chen property}. Recall  from the Introduction that $[\del f]_{tus}:=f_{ts}-f_{tu}-f_{us}$
if $f$ is a function of two variables.

\begin{Proposition}[tree Chen property] (see \cite{Kre99} or \cite{Gu2}) \label{prop:treex}

Let $\T$ be a forest, then

\BEQ [\del I_{\T}(\Gamma)]_{tus}=\sum_{\vec{v}\models V(\T)} [I_{Roo_{\vec{v}}\T}(\Gamma)]_{tu}
[I_{Lea_{\vec{v}} \T}(\Gamma)]_{us}. \label{eq:treex} \EEQ

\end{Proposition}

This proposition is illustrated in the discussion following Lemma \ref{lem:tree-skdec} in the
upcoming paragraph.

%%%%%%%%%%%%%%%%%%%%%%%%%%%%%%%%%%%%5

\subsection{Skeleton integrals}

%%%%%%%%%%%%%%%%%%%%%%%%%%%%%%%%%

We now introduce a variant of tree iterated integrals that we call {\em tree skeleton integrals}, or simply
{\em skeleton integrals}. We explain after eq. (\ref{eq:2:reg-sk-it0}) below the reason why we shall use skeleton
integrals instead of usual iterated integrals as building stones for our construction.

\begin{Definition}[formal integral]

Let $f:\R\to\R$ be a smooth, compactly supported function such that ${\cal F}f(0)=0$.
Then the formal integral $\int^t f$ of $f$ is defined as 
\BEQ \int^t f:=\frac{1}{\sqrt{2\pi}}
\int_{-\infty}^{+\infty} ({\cal F}f)(\xi) \frac{e^{\II t\xi}}{\II \xi}\ d\xi. \label{eq:intt} \EEQ

\end{Definition}

The condition ${\cal F}f(0)=0$ prevents possible infra-red divergence
when $\xi\to 0$.
Note that 
\BEQ \int^t f-\int^s f=\frac{1}{\sqrt{2\pi}} \int_{-\infty}^{+\infty} ({\cal F}f)(\xi)
\left( \int_s^t e^{\II x\xi} dx\right)d\xi=\int_s^t f(x) dx \EEQ
by the Fourier inversion formula, so $\int^t f$ is an anti-derivative of $f$.

Formally one may write, as an equality of distributions: 
\BEQ \int^t e^{\II x\xi} dx=\int_{\infty}^t e^{\II x\xi} dx=\frac{e^{\II t\xi}}{\II\xi} \EEQ
since $\int_{-\infty}^{+\infty} \frac{e^{\II x\xi}}{\II\xi} \phi(\xi)\ d\xi\to_{x\to\infty} 0$ for any test function
$\phi$ such that $\phi(0)=0$. Hence 
\BEQ \int^t f=\frac{1}{\sqrt{2\pi}} \int_{-\infty}^{+\infty} d\xi ({\cal F}f)(\xi) \int^t e^{\II x\xi}
dx=\frac{1}{\sqrt{2\pi}} \int_{-\infty}^{+\infty} ({\cal F}f)(\xi) \frac{e^{\II t\xi}}{\II\xi}
d\xi,\EEQ 
in coherence with eq. (\ref{eq:intt}).

\begin{Definition}[skeleton integrals] \label{def:skeleton-integrals}

\begin{itemize}
\item[(i)]

Let $\T$ be a tree with decoration 
$\ell:\T\to\{1,\ldots,d\}$.  Let $(v_1,\ldots,
v_{|V(\T)|})$ be any ordering of $V(\T)$ compatible with the tree partial ordering. Then the
{\em skeleton integral} of $\Gamma$ along $\T$ is by definition
\BEQ [\SkI_{\T}(\Gamma)]_t:=\int^t d\Gamma_{x_{v_1}}(\ell(v_1))
\int^{x_{v_2^-}} d\Gamma_{x_2}(\ell(v_2))
\ldots \int^{x_{v^-_{|V(\T)|}}} d\Gamma_{x_{v_{|V(\T)|}}}(\ell(v_{|V(\T)|})).
\EEQ

\item[(ii)] (extension to forests) Let $\T=\T_1\ldots \T_l$ be a forest, with its tree decomposition. Then one
defines
\BEQ [\SkI_{\T}(\Gamma)]_t:=\prod_{j=1}^l [\SkI_{\T_j}(\Gamma)]_t.\EEQ

\item[(iii)] (multilinear extension, see Definition \ref{def:2:it-int}) Assume $\T$ is
a subtree of $\tilde{\T}$, and $\mu$ a compactly supported, signed Borel measure on $\R^{\tilde{\T}}:=
\{(x_v)_{v\in V(\tilde{\T})}, x_v\in\R\}$. Then
\BEQ [\SkI_{\T}(\mu)]_t:=\int^t \int^{x_{v_2^-}}\ldots\int^{x_{v^-_{|V(\T)|}}} \mu(dx_{v_1},\ldots,
dx_{v_{|V(\T)|}}) \EEQ
is a signed Borel measure on $\{(x_v')_{v'\in V(\tilde{\T})\setminus V(\T)}, x_{v'}\in\R\}$.

\end{itemize}

\end{Definition}

Formally again, $[\SkI_{\T}(\Gamma)]_t$ may be seen as $[I_{\T}(\Gamma)]_{t,\pm\II\infty}$.
Denote by $\hat{\mu}$ the partial Fourier transform of $\mu$ with respect to $(x_v)_{v\in V(\T)}$),
 so that
\BEQ \hat{\mu}((\xi_v)_{v\in V(\T)},(dx_{v'})_{v'\in V(\T')\setminus V(\T)})=
(2\pi)^{-|V(\T)|/2} \langle \mu, \left( (x_v)_{v\in V(\T)} \mapsto e^{-\II \sum_{v\in V(\T)}
x_v \xi_v} \right) \rangle.\EEQ
Then
\BEQ [\SkI_{\T}(\mu)]_t=(2\pi)^{-|V(\T)|/2} \langle \hat{\mu}, \left[
\SkI_{\T}\left( (x_v)_{v\in V(\T)} \mapsto  e^{\II \sum_{v\in V(\T)} x_v \xi_v} \right) \right]_t \rangle.
\label{eq:2:reg-sk-it0}  \EEQ

\bigskip

As explained in the previous subsection, tree skeleton integrals are straightforward generalizations
of usual tree iterated integrals. They are very natural when computing in Fourier coordinates, because every
successive integration brings about a new $\xi$-factor in the denominator, allowing easy H\"older estimates
using Besov norms (see Appendix). On the contrary, $\int_0^t e^{\II x\xi} dx=\frac{e^{\II t\xi}}{\II \xi}-\frac{1}{\II\xi}$
contains a constant term $-\frac{1}{\II\xi}$ which does not improve when one integrates again.

\medskip

It is the purpose of section 3 to show that a rough path $\bf\Gamma$
over an $\alpha$-H\"older path $\Gamma$ may be obtained from adequately {\em regularized} tree
skeleton integrals, using the following {\em tree skeleton decomposition}, which is a variant of
 the {\em tree Chen property} recalled in Proposition \ref{prop:treex} above.

\begin{Definition}[multiple cut] \label{def:multiple-cut}

Let $\vec{v}\subset V(\T)$, $\vec{v}\not=\emptyset$. If $w\in \vec{v}$, one calls $Lev(w):=1+|\{w'\in\vec{v};
w\twoheadrightarrow w'\}|$ the {\em level} of $w$. If $\vec{v}\models V(\T)$ is an admissible cut, then $Lev(w)=1$
for all $w\in\vec{v}$. Quite generally, letting $Lev(\vec{v})=\max\{Lev(w); w\in\vec{v}\}$, one writes
$\vec{v}_j:=\{w\in\vec{v}; Lev(w)=j\}$ for $1\le j\le Lev(\vec{v})$, and calls $(\vec{v}_j)_{j=1,\ldots,Lev(\vec{v})}$
the {\em level decomposition} of $\vec{v}$ considered as a multiple cut. One shall also write: $\vec{v}_1\models
\ldots\models \vec{v}_{Lev(\vec{v})}\models V(\T)$ since $\vec{v}_{Lev(\vec{v})}\models V(\T)$ and each $\vec{v}_j$,
$j=1,\ldots,Lev(\vec{v})-1$ is an admissible cut of $Roo_{\vec{v}_{j+1}}(\T)$. 

\end{Definition}

\begin{Lemma}[tree skeleton decomposition] \label{lem:tree-skdec}

Let $\T$ be a tree. Then:

\begin{itemize}
\item[(i)] ({\em recursive version})
\BEQ [I_{\T}(\Gamma)]_{tu}=[\del \SkI_{\T}(\Gamma)]_{tu}-\sum_{\vec{v}\models V(\T)}
[I_{Roo_{\vec{v}\T}}(\Gamma)]_{tu} \ .\ [\SkI_{Lea_{\vec{v}\T}}(\Gamma)]_u, \label{eq:skdec} \EEQ

\item[(ii)]  ({\em non-recursive version})

\BEA && [I_{\T}(\Gamma)]_{tu}=[\del \SkI_{\T}(\Gamma)]_{tu}+ 
 \sum_{l\ge 1} \sum_{\vec{v}_1\models\ldots\models \vec{v}_l\models V(\T)}
 (-1)^{|\vec{v}_1|+\ldots+|\vec{v}_l|}  \nonumber\\
&& \qquad [\del \SkI_{Roo_{\vec{v}_1}(\T)} (\Gamma)]_{tu} \prod_{m=1}^{l-1}
\left[ \SkI_{Lea_{\vec{v}_m}\circ Roo_{\vec{v}_{m+1}}(\T)}\right]_u
\ [\SkI_{Lea_{\vec{v}_l}(\T)}(\Gamma)]_u. \nonumber\\ 
\EEA

\end{itemize}
\end{Lemma}

{\bf Proof.} Same as for Proposition \ref{prop:treex}.  Eq. (\ref{eq:skdec})
 may formally be seen as a particular case of the Chen property (\ref{eq:treex}) by setting
$s=\pm\II\infty$ (see previous subsection). The non-recursive version may be deduced
from the recursive version in a straightforward way. \hfill \eop

\medskip

Let us illustrate these notions in a more pedestrian way for the reader who is not accustomed
to tree integrals.
Consider for an example the trunk tree $\T_n$ with vertices $n\to n-1\to\ldots\to 1$ and
decoration $\ell:\{1,\ldots,n\}\to\{1,\ldots,d\}$, and the associated iterated integral
\BEQ [I_n^{\ell}(\Gamma)]_{ts}=
 [I_{\T_n}(\Gamma)]_{ts}=\int_s^t d\Gamma_{x_1}(\ell(1))\ldots\int_s^{x_{n-1}}d\Gamma_{x_{n}}
(\ell(n)).\EEQ

Cutting $\T_n$ at some vertex $v\in\{2,\ldots,n\}$ produces two trees, $Roo_v \T_n$ 
 and $Lea_v\T_n$, with
respective vertex subsets $\{1,\ldots,v-1\}$ and $\{v,\ldots,n\}$. 
Then the usual Chen property (ii) in the Introduction
reads
\BEQ [\del I_{\T_n}(\Gamma)]_{tus}=\sum_{v\in V(\T_n)\setminus\{1\}} [I_{Roo_v\T_n}(\Gamma)]_{tu}
[I_{Lea_v \T_n}(\Gamma)]_{us}. \label{eq:treex0} \EEQ

On the other hand, rewrite $[I_{\T_n}(\Gamma)]_{tu}$
as the sum of the {\em increment term}, which is a skeleton integral,
\BEA  [\del \SkI_{\T_n}(\Gamma)]_{tu} & =\int^t d\Gamma_{x_1}(\ell(1))\int^{x_1} d\Gamma_{x_2}(\ell(2))
\ldots\int^{x_{n-1}}d\Gamma_{x_{n}}(\ell(n))\nonumber\\
& \qquad   - \int^u d\Gamma_{x_1}(\ell(1))\int^{x_1} d\Gamma_{x_2}(\ell(2))
\ldots\int^{x_{n-1}}d\Gamma_{x_{n}}(\ell(n)) \nonumber\\   \label{eq:1.26} \EEA
and of the {\em boundary term}
\BEA &&  [I_{\T_n}(\Gamma)(\partial)]_{tu}:=  -\sum_{n_1+n_2=n} \int_u^t d\Gamma_{x_1}(\ell(1))
\ldots \int_u^{x_{n_1-1}}d\Gamma_{x_{n_1}}(\ell(n_1)) \ .\ \nonumber\\
 && \qquad  .\ \int^u d\Gamma_{x_{n_1+1}}(\ell(n_1+1))
\int^{x_{n_1+1}}d\Gamma_{x_{n_1+2}}(\ell(n_1+2)) \ldots \int^{x_{n-1}} d\Gamma_{x_{n}}(\ell(n)). \nonumber\\
\label{eq:1.27}
\EEA

The above decomposition is fairly obvious for $n=2$  and obtained by easy induction for
general $n$. One has thus obtained the {\em recursive} skeleton decomposition property for trunk trees,
\BEQ [I_{\T_n}(\Gamma)]_{tu}=[\del \SkI_{\T_n}(\Gamma)]_{tu}-\sum_{v\in V(\T_n)\setminus\{1\}}
[I_{Roo_v\T_n}(\Gamma)]_{tu} \ .\ [\SkI_{Lea_v\T_n}(\Gamma)]_u. \label{eq:2:*0} \EEQ

The {\em non-recursive} version of the skeleton decomposition property is a straightforward consequence,
and reads in this case
\BEA && [I_{\T_n}(\Gamma)]_{tu} = [\del \SkI_{\T_n}(\Gamma)]_{tu} +\sum_{l\ge 1} (-1)^l \times
\nonumber\\ &&  \qquad \qquad \times  \sum_{j_1<\ldots<j_l} [\del \SkI_{Roo_{j_1}(\T_n)}(\Gamma)]_{tu}
\prod_{m=1}^{l-1} [\SkI_{Lea_{j_m}\circ Roo_{j_{m+1}}(\T_n)}(\Gamma)]_u [\SkI_{Lea_{j_l}(\T_n)}(\Gamma)]_u, \nonumber\\
\EEA
where $Lea_{j_m}\circ Roo_{j_{m+1}}\T_n$ is the piece of $\T_n$ with subset of vertices ranging in
$\{j_m,\ldots,j_{m+1}-1\}$.

%%%%%%%%%%%%%%%%%%%%%%%%%%%%%%%%%%%%%%%%%
%%%%%%%%%%%%%%%%%%%%%%%%%%%%%%%%

\section{Regularization : the Fourier normal ordering step by step.}

%%%%%%%%%%%%%%%%%%%%%%%%%%%%%%%%%%%%%%%%%
%%%%%%%%%%%%%%%%%%%%%%%%%%%

We now come back to the original problem and assume $\Gamma$ is a $d$-dimensional $\alpha$-H\"older,
compactly supported,
{\em non-smooth} path. Then none of the previous definitions relative to iterated integrals make sense.
However, one may rewrite these as diverging series such that every term is well-defined. This
follows easily from the Besov decomposition given in the Appendix. Let us recall briefly, refering
to the Appendix for details and notations, that $\Gamma$ may be decomposed as $\sum_{k\in\Z}
D(\phi_k)\Gamma$, where $(\phi_k)_{k\in\Z}$ is a dyadic partition of unity, and
$D(\phi_k)\Gamma={\cal F}^{-1}(\phi_k \cdot {\cal F}\Gamma)$. The Fourier transform $\cal F$ has been
 introduced
at the end of the Introduction. Since $\phi_k \cdot {\cal F}\Gamma$ is a
compactly supported $C^{\infty}$ function, 
\BEQ D(\phi_k)\Gamma:x\mapsto \frac{1}{\sqrt{2\pi}}
\int_{\R} \phi_k(\xi) ({\cal F}\Gamma)(\xi) e^{\II x\xi} d\xi \EEQ
 is a $C^{\infty}$-function,
and it makes perfectly sense to integrate the $D(\phi_k)\Gamma(i)$, $k\in\Z$, $1\le i\le d$ against
each other. We suggest the following definition, where $\T\in {\cal T}$ is a fixed tree. All $\cal P$-projections
below extend to measures $\mu\in Meas(\R^{\T})$, where $\R^{\T}:=\{(x_v)_{v\in V(\T)}, x_v\in\R\}$.

\begin{Definition}[${\cal P}$-projections] \label{def:P-projections}

\begin{itemize}
\item[(i)] Let, for $\vec{k}\in \Z^{\T}:=\{(k_v)_{v\in V(\T)}, k_v\in\Z\}$,
\BEQ {\cal P}^{\{\vec{k}\}}(\Gamma):= \otimes_{v\in V(\T)} D(\phi_{k_v})\Gamma(\ell(v)),\EEQ

\item[(ii)] Similarly,
let $U\subset\Z^{\T}$. Then  
\BEQ {\cal P}^U (\Gamma):=\sum_{\vec{k}=(k_v)_{v\in V(\T)}\in U}  {\cal P}^{\{\vec{k}\}}(\Gamma).
\EEQ
\item[(iii)]
Let in particular  ${\cal P}^{+,\T}$  be 
 the ${\cal P}$-projection associated to the subset
\BEQ U=\Z^{\T}_+:=\{(k_v)_{v\in V(\T)} \in\Z^{\T}\ |\ (v \twoheadrightarrow w)\Rightarrow |k_v|\ge |k_w|\}.\EEQ
If $\T=\T_n$ is the trunk tree with $n$ vertices $\{n\to \ldots\to 1\}$ and decoration $\ell:j\mapsto j$, $j=1,\ldots,n$, see Fig. \ref{Fig1}, we shall simply write ${\cal P}^+$ instead of ${\cal P}^{+,\T_n}$. More generally, if a tree $\T$ is equipped
with a partial or total ordering $>$ compatible with its tree ordering, we let ${\cal P}^+:={\cal P}^{U_>}$ with
$U_>:=\{(k_v)_{v\in V(\T)} \in\Z^{\T}\ |\ (v>w)\Rightarrow |k_v|\ge |k_w|\}$. 

\item[(iv)]  Using the Fourier multipliers $D(\tilde{\phi}_{k_v})$ instead of $D(\phi_{k_v})$, see Definition
\ref{def:phitildek}, define similarly 
\BEQ \tilde{\cal P}^{\{\vec{k}\}}:=\frac{1}{|\Sigma_{\vec{k}}|}
 \otimes_{v\in V(\T)} D(\tilde{\phi}_{k_v})\Gamma(\ell(v)),\EEQ 
 where $\Sigma_{\vec{k}}\subset \Sigma_n$ is the subset of permutations $\tau$ such that $|k_{\tau(j)}|=|k_j|$
for every $j=1,\ldots,n$, and
\BEQ \tilde{\cal P}^+:=\sum_{\vec{k}=(k_v)_{v\in V(\T)}\in U_>}  \tilde{\cal P}^{\{\vec{k}\}}(\Gamma).\EEQ
\end{itemize}
\end{Definition}

{\bf Remark.} By construction,
${\cal P}^+ \tilde{\cal P}^+=\tilde{\cal P}^+$ if ${\cal P}^+$, $\tilde{\cal P}^+$ are associated to a
 total ordering compatible with the tree ordering of $\T$.  

\medskip

Note that ${\cal P}^U$ may be considered as a linear operator
${\cal P}^U:(B^{\alpha}_{\infty,\infty})^{\otimes\T}\to
(B^{\alpha}_{\infty,\infty})^{\otimes\T}$, where
  $(B_{\infty,\infty}^{\alpha})^{\otimes\T}$ stands for the vector space generated by the monomials
$\otimes_{v\in V(\T)} f_v$, $f_v\in B_{\infty,\infty}^{\alpha}$. It is actually a  bounded linear 
operator, as recalled in the Appendix, see Proposition \ref{prop:app:Fourier-x} and remarks after Proposition
\ref{prop:1:phi}.

\bigskip

We may now proceed to explain our {\em regularization algorithm}.

\begin{itemize}

\item {\em Step 1 (choice of regularization scheme).}
  Choose for each tree $\T\in{\cal T}$ a subset $\Z^{\T}_{reg}\subset\Z^{\T}_+$ such that the series
$\sum_{\vec{k}\in\Z^{\T}_{reg}} [\SkI_{\T}( {\cal P}^{\{\vec{k}\}}(\Gamma))]_t$ converges absolutely
for {\em any} $\alpha$-H\"older path $\Gamma$. By assumption $\Z^{\T}_{reg}=\Z$ if $|V(\T)|=1$.

\item {\em Step 2.} Let $\T$ be a forest equipped with a partial or total ordering compatible with its
tree ordering, and  $\tilde{\cal P}^+$ the corresponding projection operator. 
 For $\vec{k}\in \Z^{\T}_+$, we let the {\em projected regularized skeleton
integral} be the quantity
\BEQ [{\cal R}^{\{\vec{k}\}}\SkI_{\T}(\tilde{\cal P}^+\Gamma)]_t={\bf 1}_{\vec{k}\in\Z^{\T}_{reg}} \cdot
 [\SkI_{\T}( {\cal P}^{\{\vec{k}\}}\tilde{\cal P}^+\Gamma)]_t.\EEQ

\item {\em Step 3 (regularized projected tree integral).} For $\vec{k}\in\Z^{\T}_+$, let
$[{\cal R}^{\{\vec{k}\}}I_{\T}(\tilde{\cal P}^+\Gamma)]_{ts}$ be constructed out of projected regularized skeleton integrals
 in the following recursive way, as in Lemma \ref{lem:tree-skdec}:

\BEA &&  [{\cal R}^{\{\vec{k}\}}I_{\T}(\tilde{\cal P}^+\Gamma)]_{ts}:=[\del {\cal R}^{\{\vec{k}\}} \SkI_{\T}(
\tilde{\cal P}^+\Gamma)]_{ts}
- \nonumber\\
&& \qquad \sum_{\vec{v}\models V(\T)} [{\cal R}^{\{Roo_{\vec{v}}(\vec{k})\}} 
I_{Roo_{\vec{v}}(\T)}(\tilde{\cal P}^+\Gamma)]_{ts} [{\cal R}^{\{Lea_{\vec{v}}(\vec{k})\}} \SkI_{Lea_{\vec{v}}\T} (
\tilde{\cal P}^+\Gamma)
]_s, \nonumber\\ \label{eq:2:skdec}  \EEA
where $Roo_{\vec{v}}(\vec{k})=(k_w)_{w\in Roo_{\vec{v}}(\T)} \in \Z^{Roo_{\vec{v}}(\T)}$,
 and $Lea_{\vec{v}}(\vec{k})
=(k_w)_{w\in Lea_{\vec{v}}(\T)} \in \Z^{Lea_{\vec{v}}(\T)}.$

\item {\em Step 4 (generalization to forests).} The generalization is straightforward. Namely,
if $\T=\T_1\ldots\T_l$ is a forest, and $\vec{k}=(\vec{k}_1,\ldots,\vec{k}_l)\in\Z_+^{\T_1}\times
\ldots\times\Z_+^{\T_l}$, we let
\BEQ {\cal R}^{\{\vec{k}\}}\SkI_{\T}(\tilde{\cal P}^+\Gamma):=\prod_{j=1}^l {\cal R}^{\{\vec{k}_j\}} \SkI_{\T_j}(
\tilde{\cal P}^+\Gamma)\EEQ
and similarly
\BEQ {\cal R}^{\{\vec{k}\}} I_{\T}(\tilde{\cal P}^+\Gamma):=\prod_{j=1}^l {\cal R}^{\{\vec{k}_j\}}  I_{\T_j}(
\tilde{\cal P}^+\Gamma).\EEQ

Consider a partial or total ordering $>$ on $\T$ and denote by $\tilde{\cal P}^+$ the corresponding
projection operator. By summing over all indices\\
 $\vec{k}\in U_>$, one gets the following quantities,
\BEQ {\cal R}\SkI_{\T}( \tilde{\cal P}^+ \Gamma):= \sum_{\vec{k}\in U_>}   {\cal R}^{\{\vec{k}\}}\SkI_{\T}(
\tilde{\cal P}^+\Gamma) \label{eq:2.9} \EEQ
(see Definition \ref{def:P-projections}),   and similarly
\BEQ {\cal R}I_{\T}( \tilde{\cal P}^+ \Gamma):= \sum_{\vec{k}\in U_>}   {\cal R}^{\{\vec{k}\}}I_{\T}(
\tilde{\cal P}^+ \Gamma). \EEQ

Observe in particular, using eq. (\ref{eq:2:skdec}), and summing over indices $\vec{k}$,
 that ${\cal R}I_{\T}(\tilde{\cal P}^+\Gamma)$
decomposes naturally into the sum of an {\em increment term}, which is a regularized skeleton integral, and of
a {\em boundary term} denoted by the symbol $\partial$, namely,
\BEQ \left[ \del {\cal R}\SkI_{\T}(\tilde{\cal P}^+ \Gamma) \right]_{ts} + \left[
 {\cal R}I_{\T}(\tilde{\cal P}^+ \Gamma)(\partial)  \right]_{ts}. \label{eq:incr-bdry} \EEQ
This decomposition is a generalization of that obtained in subsection 1.4, see eq. (\ref{eq:1.26}) and
 (\ref{eq:1.27}).
Observe also that we {\em have not defined} ${\cal R}\SkI_{\T}(\Gamma)$, nor
${\cal R}I_{\T}(\Gamma)$; the regularized integration operators ${\cal R}I_{\T}$, ${\cal R}\SkI_{\T}$
only act on {\em Fourier normal ordered projections of paths} $\tilde{\cal P}^+\Gamma$.

\item{\em Final step (Fourier normal ordering).} Let $\T_n$ be a trunk tree with $n$ vertices
decorated by $\ell$, and, for each $\sigma\in\Sigma_n$, $\T^{\sigma}=\sum_{j=1}^{J_{\sigma}}
g(\sigma,j)\T^{\sigma}_j$ be the corresponding permutation graph,  as in Lemma \ref{lem:2:sigma}. Each forest
$\T^{\sigma}$ comes with a total ordering compatible with its tree ordering, which defines a projection
operator $\tilde{\cal P}^+$; we write for short $\tilde{\cal P}^{\sigma}\Gamma$ instead of 
$\tilde{\cal P}^+
(\otimes_{m=1}^n \Gamma(\ell(\sigma(m))))$. Then we let
\BEA &&  [{\cal R}{\bf\Gamma}^n(\ell(1),\ldots,\ell(n))]_{ts}:= 
 \sum_{\sigma\in\Sigma_n} \sum_{j=1}^{J_{\sigma}} g(\sigma,j) {\cal R}I_{\T^{\sigma}_j} (\tilde{\cal P}^{\sigma}\Gamma) \nonumber\\
&&  =\sum_{\sigma\in\Sigma_n} \left( 
\sum_{\vec{k}=(k_1,\ldots,k_n)\in\Z^n\ ;\ |k_{\sigma(1)}|\le \ldots\le |k_{\sigma(n)}|}
 \sum_{j=1}^{J_{\sigma}} g(\sigma,j) [{\cal R}^{\{\vec{k}\circ\sigma\}}
I_{\T_j^{\sigma}}  (\tilde{\cal P}^{\sigma}\Gamma)]_{ts} \right). \nonumber\\  \label{eq:final-step} \EEA

\end{itemize}

We shall prove in the next section that ${\cal R}{\bf\Gamma}$ satisfies the Chen (ii) and shuffle (iii)
properties of the Introduction. The H\"older property (i) will be proved in section 4 for an adequate
choice of subdomains $\Z^{\T}_{reg}$, $\T\in{\cal T}$ satisfying in particular the property required
in Step 1.

\bigskip

Some essential comments are in order.

\medskip

1. Assume that $\Gamma$ is smooth, and {\em do not regularize}, i.e., choose $\Z^{T}_{reg}=\Z^{\T}_+$.
Then eq. (\ref{eq:2:skdec}) is a recursive definition of the {\em non-regularized} projected integral
$[I_{\T}({\cal P}^{\{\vec{k}\}} \tilde{\cal P}^+\Gamma)]_{ts}$, as follows from the tree skeleton decomposition
property, see Lemma \ref{lem:tree-skdec}. Hence the right-hand side of 
formula (\ref{eq:final-step}) reads simply
\BEQ \sum_{\sigma\in\Sigma_n}
\sum_{\vec{k}=(k_1,\ldots,k_n)\in\Z^n\ ;\ |k_{\sigma(1)}|\le \ldots\le |k_{\sigma(n)}|}
 \sum_{j=1}^{J_{\sigma}} g(\sigma,j) [
I_{\T_j^{\sigma}}({\cal P}^{\{\vec{k}\}} \tilde{\cal P}^{\sigma}\Gamma)]_{ts}.  \EEQ
But this quantity is the usual iterated integral or canonical lift of $\Gamma$, \\
$ [{\bf \Gamma}^{cano,n}(\ell(1),\ldots,\ell(n)]_{ts}$, since  \BEQ \sum_{j=1}^{J_{\sigma}} g(\sigma,j)
[I_{\T_j^{\sigma}}({\cal P}^{\{\vec{k}\}}\tilde{\cal P}^{\sigma}\Gamma)]_{ts}=
[I_{\T^{\sigma}}({\cal P}^{\{\vec{k}\}}\tilde{\cal P}^{\sigma}\Gamma)]_{ts}=[I^{\ell}_n(
{\cal P}^{\{\vec{k}\}}\tilde{\cal P}^{\sigma}\Gamma)]_{ts} \EEQ
by Lemma \ref{lem:2:sigma}, and
\BEA && \sum_{\sigma\in\Sigma_n} \sum_{\vec{k}=(k_1,\ldots,k_n)\in\Z^n\ ;\
  |k_{\sigma(1)}|\le \ldots\le |k_{\sigma(n)}|}
 {\cal P}^{\{\vec{k}\}}\tilde{\cal P}^{\sigma}(\Gamma)=\sum_{\sigma\in\Sigma_n}
{\cal P}^+ \tilde{\cal P}^+ (\otimes_{m=1}^n \Gamma(\ell(\sigma(m)))) \nonumber\\
\qquad \qquad &&=  \sum_{\sigma\in\Sigma_n}
 \tilde{\cal P}^+ (\otimes_{m=1}^n \Gamma(\ell(\sigma(m))))=\Gamma \label{eq:2.15}, \EEA
see Remark after Definition \ref{def:P-projections}.

\medskip

2. Iterated integrals of order 1, $[{\cal R}{\bf\Gamma}^1(i)]_{ts}$, $1\le i\le d$, are {\em not
regularized}, namely, $[{\cal R}{\bf\Gamma}^1(i)]_{ts}=[{\bf\Gamma}^1(i)]_{ts}=\Gamma_t(i)-\Gamma_s(i)$,
because of the assumption in Step 1 which states that $\Z^{\T}_{reg}=\Z$ if $|V(\T)|=1$. Hence ${\cal R}\Gamma$ is a rough path {\em over $\Gamma$}.

\medskip

3. We propose a reformulation of this algorithm in a Hopf algebraic language in Lemma \ref{lem:3.5} below. An
equivalent algorithm is given in Definition \ref{def:alternative}. The abstract algebraic language of section 3 turns
out to be very appropriate to prove the Chen and shuffle properties.

%%%%%%%%%%%%%%%%%%%%%%%%%%%%%%%%%%%%%%%%
%%%%%%%%%%%%%%%%%%%%%%%%%%%%%%%%%

\section{Proof of the geometric and multiplicative properties}

%%%%%%%%%%%%%%%%%%%%%%%%%%%%%%%%%%%
%%%%%%%%%%%%%%%%%%%%%%%%%%%%%%%%%%%%%%%%%%

Let $\Gamma=(\Gamma(1),\ldots,\Gamma(d))$ be an $\alpha$-H\"older path.
This section is dedicated to the proof of

\begin{Theorem} \label{th:2}

Choose for each tree $\T$ a subset  $\Z^{\T}_{reg}\subset\Z^{\T}$ such that the condition of Step 1
of the construction in section 2 is satisfied, i.e. such that the regularized rough path
${\cal R}{\bf \Gamma}$ defined in section 2 is well-defined.
Then ${\cal R}{\bf \Gamma}$  satisfies the
Chen (ii) and shuffle (iii) properties of the Introduction.

\end{Theorem}

This theorem is in fact a consequence of the following very general construction, whose essence
is really algebraic. Two Hopf algebras are involved in it: the {\em Hopf algebra of decorated rooted trees}
$\bf H$, and the {\em shuffle algebra} $\bf Sh$.  As we shall presently see, the first one is related
to the {\em Chen property}, while the second one is related to the {\em shuffle property}. The first
paragraph below is devoted to an elementary presentation of these Hopf algebras in connection with
the Chen/shuffle property. Theorem \ref{th:2} is proved in the second paragraph.

%%%%%%%%%%%%%%%%%%%%%%%%%%%%%%%%%%%%

\subsection{Hopf algebras and the Chen and shuffle properties}

%%%%%%%%%%%%%%%%%%%%%%%%%%%%%%%%%%%%

\begin{enumerate}
\item Let us first consider the Hopf algebra of decorated rooted trees, $\bf H$.
Recall  the definition of the coproduct on $\bf H$,
\BEQ \Del(\T)=e\otimes\T+\T\otimes e+\sum_{\vec{v}\models V(\T)} Roo_{\vec{v}}\T \otimes Lea_{\vec{v}}\T.\EEQ
The usual convention \cite{ConKre98,ConKre00} is to write $c$ (cut) for $\vec{v}$, $R^c (\T)$ ({\em root} part) for
$Roo_{\vec{v}}\T$, $P^c (\T)$ for $Lea_{\vec{v}}\T$ ({\em leaves}), and to reverse the order of the factors in the tensor
product.

The {\em convolution} of two linear forms $f,g$ on $\bf H$ writes
\BEQ (f\ast g)(\T)=f(\T)g(e)+f(e)g(\T)+\sum_{\vec{v}\models V(\T)} f(Roo_{\vec{v}}\T) g(Lea_{\vec{v}}\T), \quad
\T\in {\bf H}. \label{eq:H-convol} \EEQ
This notion is particularly interesting for {\em characters}. A character of $\bf  H$ is a linear map
such that $\chi(\T_1.\T_2)=\chi(\T_1).\chi(\T_2)$. If $\chi_1,\chi_2$ are two characters of $\bf  H$,
then $\chi_1\ast\chi_2$ is also a character of $\bf  H$. 

The {\em tree Chen property}, see Proposition \ref{prop:treex}, may then be stated as follows. Let $\Gamma=(\Gamma(1),
\ldots,\Gamma(d))$ be a smooth path, and 
\BEQ {\bf H}^d:=\{\T\in {\bf H}; \ell: V(\T)\to \{1,\ldots,d\}\} \EEQ
be the subspace of $\bf H$ generated by forests with decoration valued in $\{1,\ldots,d\}$.
 Now, define $I^{ts}_{\Gamma}: {\bf H}^d\to\R$ to be
the following character of $\bf  H$ (see Definition \ref{def:2:it-int})
\BEQ I^{ts}_{\Gamma}(\T)=[I_{\T}(\Gamma)]_{ts}. \EEQ
 Then (as remarked in \cite{Kre99})
\BEQ I^{ts}_{\Gamma}=I^{tu}_{\Gamma}\ast I^{us}_{\Gamma}. \label{eq:3.5} \EEQ

Generalizing this property to the multilinear setting, one may also write
\BEQ I^{ts}_{\mu}(\T)=(I^{tu}\ast I^{us})_{\mu}(\T):=I^{tu}_{\mu}(\T)+I^{us}_{\mu}(\T)+\sum_{\vec{v}\models
V(\T)} I^{tu}_{Roo_{\vec{v}}(\mu)} (Roo_{\vec{v}}(\T)) I^{us}_{Lea_{\vec{v}}(\mu)} (Lea_{\vec{v}}(\T)) \EEQ
for a tensor measure $\mu=\otimes_{v\in V(\T)} \mu_v$, where $Roo_{\vec{v}}(\mu):=\otimes_{v\in V(Roo_{\vec{v}}(\T))} \mu_v$, $Lea_{\vec{v}}(\mu):=\otimes_{v\in V(Lea_{\vec{v}}(\T))}\mu_v$, and
\BEQ I^{ts}_{\mu}(\T):=(I^{tu}\ast I^{us})_{\mu}(\T):=\sum_{\vec{k}} (I^{tu}\ast I^{us})_{\mu_{\vec{k}}}(\T)
\label{eq:3.7} \EEQ
for a more general measure $\mu:=\sum_{\vec{k}}\mu_{\vec{k}}$, where each $\mu_{\vec{k}}$ is a tensor
measure. Later on we shall use these formulas for $\mu_{\vec{k}}={\bf 1}_{\vec{k}\in\Z^{\T}_+} 
d{\cal P}^{\{\vec{k}\}}(\Gamma)$ or ${\bf 1}_{\vec{k}\in\Z^{\T}_{reg}} 
d{\cal P}^{\{\vec{k}\}}(\Gamma).$

\medskip

  As for the {\em antipode} $S$, it is the multiplicative morphism $S:{\bf H}\to{\bf H}$ defined
inductively on tree generators $\T$ by (see \cite{ConKre98}, p. 219)
\BEQ S(e)=e; \quad S(\T)=-\T-\sum_{\vec{v}\models V(\T)} Roo_{\vec{v}}\T. S(Lea_{\vec{v}}\T). \label{eq:3.8} \EEQ
Applying iteratively the second relation yields an expression of $S(\T)$ in terms of multiple cuts
of $\T$ obtained by 'chopping' it \cite{ConKre98}, see Definition \ref{def:multiple-cut}, namely,
\BEA &&  S(\T)=-\T-\sum_{l\ge 1} \sum_{\vec{v}_1\models\ldots\models\vec{v}_l\models V(\T)}
(-1)^{|\vec{v}_1|+\ldots+|\vec{v}_l|}  \nonumber\\
&& \qquad Roo_{\vec{v}_1}(\T)  \left\{ \prod_{m=1}^{l-1}
Lea_{\vec{v}_m}\circ Roo_{\vec{v}_{m+1}}(\T) \right\}
Lea_{\vec{v}_l}(\T). \nonumber\\  \label{eq:3:S}
\EEA

Let $\chi_1,\chi_2$ be two characters of $\bf H$. Recall that $\chi_2\circ S$ is  the {\em convolution
inverse} of $\chi_2$,
namely, $\chi_2\circ S$ is a character and $\chi_2\ast(\chi_2\circ S)=\bar{e}$, where $\bar{e}$ is the
counity of $\bf H$, defined on generators by  $\bar{e}(e)=1$ and $\bar{e}(\T)=0$ if $\T$ is a forest.
Now eq. (\ref{eq:H-convol}) and  (\ref{eq:3:S}) yield
\BEA && \chi_1\ast (\chi_2\circ S)(\T)=\chi_1(\T)+\chi_2\circ S(\T)+\sum_{\vec{v}\models V(\T)}
\chi_1(Roo_{\vec{v}}(\T))\chi_2\circ S(Lea_{\vec{v}}(\T)) \nonumber\\
&&=(\chi_1-\chi_2)(\T)+\sum_{\vec{v}\models V(\T)} (\chi_1-\chi_2)(Roo_{\vec{v}}(\T))\chi_2\circ S(
Lea_{\vec{v}}(\T)) \nonumber\\
&&=(\chi_1-\chi_2)(\T)+\sum_{l\ge 1} (-1)^{|\vec{v}_1|+\ldots+|\vec{v}_l|} \sum_{\vec{v}=(\vec{v}_1,\ldots,\vec{v}_l)} 
(\chi_1-\chi_2)(Roo_{\vec{v}_1}(\T)) \times \nonumber\\
&& \qquad \qquad \qquad \qquad \qquad  \times \left[\prod_{m=1}^{l-1} \chi_2(Lea_{\vec{v}_m}\circ Roo_{\vec{v}_{m+1}(\T)})
\right]
\chi_2(Lea_{\vec{v}_l}(\T)) \nonumber\\
\EEA
where $\vec{v}=(\vec{v}_1,\ldots,\vec{v}_l)$ is a multiple cut of $\T$ as in eq. (\ref{eq:3:S}).

In particular, let $\SkI_{\Gamma}^t:{\bf H}\to\R$ be the character defined by (see Definition \ref{def:skeleton-integrals})
\BEQ \SkI_{\Gamma}^t(\T)=[\SkI_{\T}(\Gamma)]_t.\EEQ
Then the tree skeleton decomposition, see Lemma \ref{lem:tree-skdec}, reads simply
\BEQ I_{\Gamma}^{tu}=\SkI_{\Gamma}^t\ast \left(\SkI_{\Gamma}^u\circ S \right). \label{eq:3.12} \EEQ

\item

The {\em shuffle algebra}  over the index set $\N$ \cite{Mur}  may be
defined as follows. The algebra
$\bf Sh$ is generated as a vector space over $\R$
 by the identity $e$ and by the  trunk trees $(\T_n)_{n\ge 1}$ with vertex set
$V(\T_n)=\{v_1<\ldots<v_n\}$, provided with  an $\N$-valued decoration $\ell$.  Let $\T_n$,
$\T'_{n'}$ be trunk trees with $n$, resp. $n'$ vertices. The {\em shuffle
product} of $\T_n$ and $\T'_{n'}$ is the formal sum
\BEQ \T_n \shuffle \T'_{n'}=\sum_{\eps\in Sh((V(\T_n),V(\T'_{n'})))} \eps(\T_n^{\T'_{n'}}),\EEQ
where $\T_n^{\T'_{n'}}$ is the trunk tree with $n+n'$ vertices obtained by putting $\T'_{n'}$ on top
of $\T_n$, and the shuffle $\eps$ permutes the decorations of $\T_n$, $\T'_{n'}$ as in property
(iii) discussed in the Introduction.

Let ${\bf Sh}^d$ be the subspace of $\bf Sh$ generated by trunk trees with decoration valued in $\{1,\ldots,d\}$.
 Then the {\em shuffle property} for iterated integrals
reads
\BEQ I_{\Gamma}^{ts}(\T_n) I_{\Gamma}^{ts}(\T'_{n'}) = I_{\Gamma}^{ts}(\T_n\shuffle\T'_{n'}), \quad
\T_n,\T'_{n'}\in {\bf Sh}^d. \label{eq:3.14} \EEQ
In other words, it may be stated by saying that $I_{\Gamma}^{ts}:\T_n\to [I_{\T_n}(\Gamma)]_{ts}$ is
{\em a character of} $\bf Sh$. Similarly, {\em skeleton integrals} $\SkI_{\Gamma}^t:\T_n\to [\SkI_{\T}(\Gamma)]_t$ also
define characters of $\bf Sh$.

The shuffle algebra $\bf Sh$ is made into a Hopf algebra by re-using the same coproduct
$\Del:\T\to \T\otimes e+e\otimes \T+\sum_{\vec{v}\models V(\T)} Roo_{\vec{v}}\T\otimes Lea_{\vec{v}}\T$ as
for $\bf H$, and defining the {\em antipode} $\bar{S}$
 as $\bar{S}(\T_n)=(-1)^n \bar{\T}_n$, where $\bar{\T}_n$
is obtained from  $\T_n$ by reversing the  ordering of the vertices,
 $\ell_{\bar{\T}_n}(v_j)=\ell_{\T_n}(v_{n+1-j})$.

The convolution of linear forms or characters $f,g$ on ${\bf Sh}$ is given by the same formula
as for $\bf H$.

\begin{Proposition} \cite{Mur}  \label{prop:Pi}

The linear morphism $\Pi:{\bf H}\to {\bf Sh}$ defined by $ \Pi(\T)=\sum_j \T_j,$ where $\T_j$ ranges
over all trunk trees $\{v_1<\ldots<v_{|V(\T)|}\}$ such that the corresponding
total ordering of vertices of $\T$
is compatible with its tree partial ordering, is a Hopf algebra map.

\end{Proposition}

$\Pi$ is actually {\em onto}. In other words, it is a {\em structure-preserving projection},
 with the canonical identification of $\bf Sh$
as a subspace of $\T$. Note that $[I_{\T}(\Gamma)]_{ts}=[\SkI_{\T}(\Gamma)]_{ts}=0$ if $\T\in Ker(\Pi)$ and
$\Gamma$ is an arbitrary smooth path, which is a straightforward generalization of the shuffle property; one
may call this the {\em tree shuffle property}.

\begin{Corollary} \label{cor:chiS-chibarSbar}

Let $\bar{\chi}$ be a caracter of $\bf Sh$. Then $\chi:=\bar{\chi}\circ\Pi$ is a character of
$\bf H$. If $\T\in {\bf Sh}$, then $\chi\circ S(\T)=\bar{\chi}\circ\bar{S}(\T)$.

\end{Corollary}

\end{enumerate}

%%%%%%%%%%%%%%%%%%%%%%%%%%%%%%%%%%

\subsection{Proof of the Chen and shuffle properties}

%%%%%%%%%%%%%%%%%%%%%%%%%%%%%%%

We shall now prove Theorem \ref{th:2}. In the next pages, $Meas(\R^n)$ stands for the space of compactly
supported, signed Borel measures on $\R^n$. Let us explain the strategy of the  proof. We give a general method
to construct families of characters of the shuffle algebra, $\bar{\chi}^t_{\Gamma}$, depending on a path $\Gamma$,
see Lemma \ref{lem:barchi-char}; these quantities satisfy the shuffle property by eq. (\ref{eq:3.14}).  Then
$\bar{\chi}^t_{\Gamma}\ast(\bar{\chi}^s\circ\bar{S})$ is immediately seen to define a rough path satisfying both
the Chen and shuffle properties, see Definition \ref{def:alternative}. For a particular choice of the characters
$\bar{\chi}^t_{\Gamma}$ related to the regularized skeleton integrals defined in section 2, the rough path of
Definition \ref{def:alternative} is shown to coincide with the regularized rough path ${\cal R}{\bf\Gamma}$ of 
section 2, see Lemma \ref{lem:R'=R}. In order to prove this last lemma, one needs a Hopf algebraic reformulation
of the Fourier normal ordering algorithm leading to ${\cal R}{\bf\Gamma}$, see Lemma \ref{lem:3.5}.

\medskip

\begin{Lemma}[measure splitting] \label{lem:measure-splitting}

Let $\mu\in Meas(\R^n)$. Then \BEQ \mu=\sum_{\sigma\in\Sigma_n} \mu^{\sigma}\circ\sigma,\EEQ
where $\mu^{\sigma}\in \tilde{\cal P}^+ Meas(\R^n)$ is defined by
\BEQ \mu^{\sigma}:=\sum_{\vec{k}=(k_1,\ldots,k_n)\in\Z^n; |k_{\sigma(1)}|\le \ldots\le |k_{\sigma(n)}|}
 (\tilde{\cal P}^{\{\vec{k}\}}\mu)\circ\sigma \EEQ
as in eq. (\ref{eq:final-step}).

\end{Lemma}

{\bf Proof.} See eq.  (\ref{eq:2.15}). \hfill \eop

\begin{Definition} \label{def:barchi-char}

\begin{itemize}
\item[(i)] Let ${\cal F}^+_{n,n}\subset {\bf H}$ $(n\ge 1)$ be the set of all forests
$\T$ with $n$ vertices and one-to-one decoration $\ell:V(\T)\to\{1,\ldots,n\}$ valued in
the set $\{1,\ldots,n\}$, such that $(v\twoheadrightarrow w)\Rightarrow \ell(v)\ge \ell(w)$, and 
${\bf H}^+_{n,n}\subset{\bf H}$ the vector space generated by ${\cal F}^+_{n,n}$.
\item[(ii)] If $\T\in {\cal F}^+_{n,n}$, let $\tilde{\cal P}^{+,\T} Meas(\R^n)$ denote the subspace
 $\{ \tilde{\cal P}^{+,\T}\mu; \mu\in Meas(\R^n)\}$, see section 2 for a definition of the
 projection operator
$\tilde{\cal P}^{+,\T}$.

\item[(iii)] Let $\phi^t_{\T}: \tilde{\cal P}^{+,\T}Meas(\R^n)\to\R, \mu\mapsto \phi^t_{\T}(\mu)$, also written $\phi^t_{\mu}(\T)$ 
$(t\in\R, \T\in {\cal F}^+_{n,n})$  be a family of linear forms such that, if $(\T_i,\mu_i)\in
{\cal F}^+_{n_i,n_i}\times \tilde{\cal P}^{+,\T_i} Meas(\R^{n_i})$, $i=1,2$, the following 
{\em ${\bf H}$-multiplicative property} holds,
\BEQ \phi^t_{\mu_1}(\T_1)  \phi^t_{\mu_2}(\T_2)=\phi^t_{\mu_1\otimes\mu_2}(\T_1\wedge\T_2),
\label{eq:T1wedgeT2}\EEQ
where $\T_1\wedge \T_2\in {\cal F}^+_{n_1+n_2,n_1+n_2}$ is the forest $\T_1.\T_2$ with decoration
$\ell\big|_{\T_1}=\ell_1$ , $\ell\big|_{\T_2}=n_1+\ell_2$  ($\ell_i$ = decoration of $\T_i$, $i=1,2$), and
$\mu_1\otimes\mu_2\in \tilde{\cal P}^{+,\T_1\wedge\T_2} Meas(\R^{n_1+n_2})$ is the tensor measure
$$\mu_1\otimes \mu_2(dx_1,\ldots,dx_{n_1+n_2})=\mu_1(dx_1,\ldots,dx_{n_1})\mu_2(dx_{n_1+1},\ldots,
dx_{n_1+n_2}).$$

\item[(iv)] Let, for $\Gamma=(\Gamma(1),\ldots,\Gamma(d))$,
 $\bar{\chi}_{\Gamma}^t:{\bf Sh}^d\to\R$ be the linear form on ${\bf Sh}^d$ defined by
\BEQ \bar{\chi}_{\Gamma}^t (\T_n):=\sum_{\sigma\in \Sigma_n} \phi^t_{\mu_{\Gamma}^{\sigma}}(\T^{\sigma}),\EEQ
where -- $\ell$ being the decoration of $\T_n$ -- one has set  $\mu_{\Gamma}:=\otimes_{j=1}^n d\Gamma(\ell(j))$, and $\T^{\sigma}$ is the permutation graph associated to $\sigma$ (see subsection 1.2). 
\end{itemize}

\end{Definition}

{\bf Remarks.}

\begin{enumerate}
\item Note that the ${\bf H}$-multiplicative property (\ref{eq:T1wedgeT2}) holds in particular for 
$\phi^t_{\T}
=[\SkI_{\T}(\ .\ )]_t$ or $[{\cal R}\SkI_{\T}(\ .\ )]_t$, either trivially or by construction
 (see Step 4 in the construction of section 2).   Note that $[{\cal R}\SkI_{\T}(\mu)]_t$ has been defined 
{\em only} if
$\mu\in \tilde{\cal P}^+ Meas(\R^n)$.
  If $\phi^t_{\T}=
[\SkI_{\T}(\ .\ )]_t$, then simply $\bar{\chi}^t_{\Gamma}(\T_n)=[\SkI_{\T_n}(\Gamma)]_t$
 by the measure splitting lemma.

\item Assume  $\mu_i\in \tilde{\cal P}^+ Meas(\R^{n_i})\subset \tilde{\cal P}^{+,\T}Meas(\R^{n_i})$, where
$\tilde{\cal P}^+$ is the $\tilde{\cal P}$-projection associated to the subset $\Z^{n_i}_+:=\{\vec{k}=(k_1,\ldots,k_{n_i});
|k_1|\le \ldots\le |k_{n_i}|\}$ $(i=1,2)$. Then $\mu_1\otimes\mu_2\in \tilde{\cal P}^{+,\T_1\wedge\T_2}Meas(\R^{n_1+n_2})$
but $\mu_1\otimes\mu_2\not\in \tilde{\cal P}^+ Meas(\R^{n_1+n_2})$ in general; the product measure $\mu_1\otimes\mu_2$
decomposes as a sum over shuffles $\eps$ of $(1,\ldots,n_1),(n_1+1,\ldots,n_1+n_2)$, namely,
$\mu_1\otimes\mu_2=\sum_{\eps\ {\mathrm{shuffle}}} (\mu_1\otimes\mu_2)^{\eps}\circ\eps$. Hence the
${\bf H}$-multiplicative property (\ref{eq:T1wedgeT2}) reads also
\BEQ  \phi^t_{\mu_1}(\T_1) \phi^t_{\mu_2}(\T_2)=\sum_{\eps\ {\mathrm{shuffle}}}
\phi^t_{(\mu_1\otimes\mu_2)^{\eps}}(\eps^{-1}(\T_1\wedge\T_2)),
\label{eq:T1wedgeT2bis}\EEQ
where $\eps^{-1}(\T_1\wedge\T_2)$ is the forest $\T_1\wedge\T_2$ with decoration $\eps^{-1}\circ\ell$, see
Definition \ref{def:barchi-char} (iii) for the definition of $\ell$.

\item The regularization algorithm $\cal R$ presented in section 2 may we written in a compact way
using the structures we have just introduced. Namely, one has:

\begin{Lemma} \label{lem:3.5}

Let $\Gamma=(\Gamma(1),\ldots,\Gamma(d))$ and $\mu_{\Gamma}:=\otimes_{j=1}^n d\Gamma(\ell(j))$. Then
\BEQ [{\cal R}{\bf\Gamma}^n(\ell(1),\ldots,\ell(n))]_{ts}=\sum_{\sigma\in\Sigma_n}
\left( \phi^t\ast(\phi^s\circ S)\right)_{\mu_{\Gamma}^{\sigma}}(\T^{\sigma}),\EEQ
where 
\BEQ \phi^t_{\nu}(\T):=[{\cal R}\SkI_{\T}(\nu)]_t=\left[\SkI_{\T} \left( 
  \sum_{\vec{k}\in\Z^{\T}_{reg}} \left( \otimes_{v\in V(\T)} D(\phi_{k_v})\right) \nu \right) \right]_t \EEQ
 for $\nu\in \tilde{\cal P}^{+,\T}Meas(\R^n)$, and
$\left(\phi^t\ast (\phi^s\circ S)\right)_{\mu^{\sigma}}$ is the obvious multilinear extension of
the convolution, see  eq. (\ref{eq:3.7}).

\end{Lemma}

{\bf Proof.} Simple formalization of the regularization procedure explained in Section 2. \hfill \eop

\end{enumerate}

\medskip

The fundamental result is the following.

\begin{Lemma} \label{lem:barchi-char}
Let $\Gamma=(\Gamma(1),\ldots,\Gamma(d))$ be compactly supported, and assume that the condition of Step 1
in section 2 is satisfied. Then  $\bar{\chi}_{\Gamma}^t$ is a character of ${\bf Sh}^d$.
\end{Lemma}

{\bf Proof.}  Let $\T_{n_i}\in {\bf Sh}^d$ 
with $n_i$ vertices $(i=1,2)$; define $n:=n_1+n_2$.  Let $\mu_i:=\otimes_{j=1}^{n_i} d\Gamma(\ell_i(j))$,
$i=1,2$ and $\mu:=\mu_1\otimes\mu_2$. If $n'\ge 1$, we let  $\T'_{n'}$  be the trunk tree with $n'$ vertices
$\{n'\to\ldots \to 1\}$ and decoration $\ell(j)=j$, $j\le n'$, see Fig. \ref{Fig1}. All shuffles $\eps$ below are
intended to be shuffles of $(1,\ldots,n_1),(n_1+1,\ldots,n_2)$.   Then
\BEA \bar{\chi}^t_{\Gamma}(\T_{n_1}\shuffle \T_{n_2}) &=& \sum_{\eps \ {\mathrm{shuffle}}}
\bar{\chi}^t_{\mu\circ \eps} (\T'_n) \nonumber\\
& = & \sum_{\sigma\in \Sigma_n} \sum_{\eps \ {\mathrm{shuffle}}} \phi^t_{(\mu\circ\eps)^{\sigma}}
(\T^{\sigma}) =\sum_{\sigma\in\Sigma_n} \sum_{\eps \ {\mathrm{shuffle}}}  \phi^t_{\mu^{\eps\circ\sigma}} (\T^{\sigma})
\nonumber\\
&=: &  \sum_{\sigma\in\Sigma_n} \phi^t_{\mu^{\sigma}}(\tgoth_1^{\sigma})  \label{eq:proof-1}  \EEA
with \BEQ \tgoth_1^{\sigma}:=\sum_{\eps \ {\mathrm{shuffle}}} \T^{\eps^{-1}\circ\sigma} \in {\bf H}^+_{n,n}.\EEQ

On the other hand,
\BEA && \bar{\chi}^t_{\Gamma}(\T_{n_1})\bar{\chi}^t_{\Gamma}(\T_{n_2})= \bar{\chi}^t_{\mu_1}(\T'_{n_1})
\bar{\chi}_{\mu_2}^t(\T'_{n_2}) \nonumber\\
&&= \sum_{\sigma_1\in\Sigma_{n_1},\sigma_2\in\Sigma_{n_2}} \phi^t_{\mu_1^{\sigma_1}}(\T^{\sigma_1})
\phi^t_{\mu_2^{\sigma_2}}(\T^{\sigma_2}) \nonumber\\
&&=\sum_{\sigma_1\in\Sigma_{n_1},\sigma_2\in\Sigma_{n_2}} \sum_{\eps\ {\mathrm{shuffle}}} \phi^t_{(\mu_1^{\sigma_1}
\otimes\mu_2^{\sigma_2})^{\eps}}(\eps^{-1}(\T^{\sigma_1}\wedge\T^{\sigma_2})) \EEA
by (\ref{eq:T1wedgeT2bis})
\BEQ =\sum_{\sigma\in\Sigma_n} \phi^t_{\mu^{\sigma}}(\tgoth_2^{\sigma}) \label{eq:proof-2}  \EEQ
where
\BEQ \tgoth_2^{\sigma}:=\sum_{(\sigma_1,\sigma_2,\eps); (\sigma_1\otimes\sigma_2)\circ\eps=\sigma}
\eps^{-1}(\T^{\sigma_1}\wedge\T^{\sigma_2}).\EEQ

Hence $\bar{\chi}^t_{\Gamma}$ is a character of ${\bf Sh}$ if and only if $\tgoth_1^{\sigma}=\tgoth_2^{\sigma}$ for every
$\sigma\in\Sigma_n$; let us prove this. Extend first (\ref{eq:proof-1}) and (\ref{eq:proof-2}) by multilinearity
from tensor measures $\mu_1\otimes \mu_2$ to a general measure $\mu\in Meas(\R^n)$. By the
usual shuffle identity, $\SkI^t_{\Gamma}(\T_{n_1}\shuffle\T_{n_2})=\SkI^t_{\Gamma}(\T_{n_1}).\SkI^t_{\Gamma}(\T_{n_2})$, 
so (\ref{eq:proof-1}) and (\ref{eq:proof-2}) coincide for $\bar{\chi}^t=[ \SkI(\ .\ )]_t$.
Choose $\sigma\in\Sigma_n$. For any $\mu\in Meas(\R^n)$, one has
\BEQ [\SkI_{\mu^{\sigma}}(\tgoth_1^{\sigma}-\tgoth_2^{\sigma})]_t=0. \EEQ
 This fact implies actually that $\tgoth_1^{\sigma}=\tgoth_2^{\sigma}$.
Let us first give an informal proof of this statement. To begin with,
 note that the fact that $[\SkI_{\Gamma}(\tgoth)]_t=0$
for every smooth path $\Gamma$ does not imply in itself that $\tgoth=0$ if $\tgoth\in {\bf H}$ is arbitrary.
Namely, the character $\SkI_{\Gamma}^t:{\bf H}\to\R$ quotients out via the canonical projection $\Pi:{\bf H}\to{\bf Sh}$, see Proposition \ref{prop:Pi},
into a character ${\bf Sh}\to\R$, by the tree shuffle property;
 one may actually prove that $\SkI^t_{\Gamma}(\tgoth)=0$ for every smooth
path $\Gamma$ if and only
if $\tgoth\in Ker(\Pi)$. In our case, the elements of ${\cal F}_{n,n}^+$ are linearly independent modulo
 $Ker(\Pi)$ because
the ordering of the labels $\ell(j)$, $j=1,\ldots,n$ is compatible with the tree ordering -- which
{\em prevents any possibility of shuffling} --, hence $\tgoth_1^{\sigma}-\tgoth_2^{\sigma}=0$.

\medskip

Let us now give a more formal argument. Let $\tgoth_1^{\sigma}-\tgoth_2^{\sigma}=:\sum_j a_j \tgoth_j$, $a_j\in\Z$, 
$\tgoth_j\in {\cal F}^+_{n,n}$ two-by-two distinct,
 and define 
\BEQ F_{\tgoth}(\xi_1,\ldots,\xi_n):=\frac{1}{\prod_{v\in V(\tgoth)} (\xi_v+\sum_{w\twoheadrightarrow v} \xi_w) }
\EEQ  if $\tgoth\in {\cal F}_{n,n}^+$. Applying Lemma \ref{lem:2:SkI} to $[\SkI_{\mu_m}(\tgoth_j)]_t$
where $(\mu_m\circ\sigma)_{m\ge 1}\in {\cal P}^+ Meas(\R^n)$ is a sequence of measures whose Fourier
transform converges weakly to the Dirac distribution $\del_{(\xi_1,\ldots,\xi_n)}$, one gets
\BEQ \sum_J a_j F_{\tgoth_j}(\xi_1,\ldots,\xi_n)=0, \quad |\xi_1|\le \ldots\le|\xi_n|. \label{eq:Ftj} \EEQ
Since the left-hand side of (\ref{eq:Ftj}) is a rational function, the equation extends to arbitrary $\vec{\xi}=(\xi_1,
\ldots,\xi_n)\in\R^n$. Note that 
\BEQ \prod_{v\in V(\tgoth_j)} (\xi_v+\sum_{w\twoheadrightarrow v} \xi_w) =(\xi_1+\sum_{w\twoheadrightarrow 1} \xi_w)
F_{\check{\tgoth}_j}(\xi_2,\ldots,\xi_n),\EEQ
where $\check{\tgoth}_j:=Lea_{\{1\}}(\tgoth_j)$ is $\tgoth_j$ severed of the vertex 1, which is 
one of its roots. Let $J_{\Omega}$, $\Omega\subset\{2,\ldots,n\}$
be the subset of indices $j$ such that $\{v\in\{1,\ldots,n\};v\twoheadrightarrow 1\ {\mathrm{in}}\ \tgoth_j\}=\Omega$, i.e.
such that the tree component of $1$ in $\tgoth_j$ has vertex set $\Omega$. Take the residue at $-\sum_{w\in\Omega}\xi_w$
of the left-hand side of (\ref{eq:Ftj}), considered as a function of $\xi_1$. This gives:
\BEQ \sum_{j\in J_{\Omega}} a_j F_{\check{\tgoth}_j}(\xi_2,\ldots,\xi_n)=0,\quad \Omega\subset\{2,\ldots,n\}.\EEQ
Shifting by $-1$ the indices of vertices of $\check{\tgoth}_j$ and the labels $\ell(v), v\in V(\check{\tgoth}_j)$, one gets a forest in ${\cal F}^+_{n-1,n-1}$.
One may now conclude by an inductive argument.

 \hfill \eop

\medskip

Let us now give an alternative definition
for the regularization  ${\cal R}$. As we shall see in Lemma \ref{lem:R'=R}, the two definitions actually coincide.

\begin{Definition}[alternative definition for regularization ${\cal R}'$] \label{def:alternative}

Choose for every tree $\T\in{\bf H}$ a subset $\Z^{\T}_{reg}\subset\Z^{\T}_+$ satisfying the condition stated in
Step 1 of section 2.  Let $\Gamma=(\Gamma(1),\ldots,\Gamma(d))$ be a compactly supported, $\alpha$-H\"older
path, and $\mu_{\Gamma}:=\otimes_{j=1}^n d\Gamma(\ell(j))$ the corresponding measure.

\begin{itemize}
\item[(i)] Let, for every $\T\in {\bf H}^d$ with $n$ vertices,
 \BEQ \phi^t_{\nu}(\T)=[{\cal R}\SkI_{\T}(\nu)]_t, \quad \nu\in\tilde{\cal P}^{+,\T} Meas(\R^n)  \EEQ
 see eq. (\ref{eq:2.9}) or Lemma \ref{lem:3.5}, and
\BEQ \bar{\chi}^t_{\Gamma}(\T_n):=\sum_{\sigma\in \Sigma_n} \phi^t_{\mu_{\Gamma}^{\sigma}}(\T^{\sigma}) \EEQ
be the associated character of $\bf Sh$ as in Definition \ref{def:barchi-char}.
\item[(ii)] Let, for $\T_n\in {\bf Sh}^d$,  $n\ge 1$, with $n$ vertices and decoration $\ell$,
  \BEQ [{\cal R}'{\bf\Gamma}^n(\ell(1),\ldots,\ell(n))]_{ts}:=
   \bar{\chi}^t_{\Gamma} \ast (\bar{\chi}^s_{\Gamma}\circ\bar{S})
(\T_n). \EEQ
\end{itemize}

\end{Definition}

Since $\bar{\chi}^s_{\Gamma}$, $ \bar{\chi}^t_{\Gamma}$ and hence $\bar{\chi}_{\Gamma}^t\ast (\bar{\chi}^s_{\Gamma}
\circ \bar{S})$  are characters of the shuffle algebra, ${\cal R}' {\bf
\Gamma}$ satisfies
the {\em shuffle property}. Also, ${\cal R'}{\bf\Gamma}$ satisfies the {\em Chen  property}  by construction,
since 
\BEA &&  [{\cal R}' {\bf\Gamma}^n(\ell(1),\ldots,\ell(n))]_{ts} =  \left(  \bar{\chi}^t_{\Gamma} \ast (\bar{\chi}^u_{\Gamma}\circ\bar{S}) \right) \ast
\left(  \bar{\chi}^u_{\Gamma} \ast (\bar{\chi}^s_{\Gamma} \circ\bar{S}) \right) (\T_n) \nonumber\\
&& \qquad = [{\cal R}'{\bf\Gamma}^n(\ell(1),\ldots,\ell(n))]_{tu}+[{\cal R}'{\bf\Gamma}^n(\ell(1),\ldots,
\ell(n))]_{us} \nonumber\\
&& \qquad \qquad \qquad \qquad +\sum_j [{\cal R}'{\bf\Gamma}^j(\ell(1),\ldots,\ell(j))]_{tu} [{\cal R}'
 {\bf\Gamma}^{n-j}(\ell(j+1),\ldots,\ell(n))]_{us} \nonumber\\ \EEA
by definition of the convolution in $\bf Sh$. Both properties remain valid if $\bar{\chi}^t_{\Gamma}$, $t\in\R$
are arbitrary characters of $\bf Sh$.

\bigskip

Let us make this definition a little more explicit before proving that ${\cal R}'={\cal R}$. 
Replacing $\bar{\chi}^s\circ\bar{S}$ with $\chi^s \circ S$, see 
Corollary \ref{cor:chiS-chibarSbar}, one gets, see
eq. (\ref{eq:3.8}),  
\BEA [{\cal R}' {\bf\Gamma}^n(\ell(1),\ldots,\ell(n))]_{ts} &=& \chi^t_{\Gamma}(\T_n)+\chi^s_{\Gamma}(S(\T_n))+\sum_j
\chi^t_{\Gamma}(Roo_j\T_n) (\chi^s_{\Gamma}\circ S)(Lea_j\T_n) \nonumber\\
&=& (\bar{\chi}^t_{\Gamma}-\bar{\chi}^s_{\Gamma})(\T_n)+\sum_j (\bar{\chi}^t_{\Gamma}-\bar{\chi}^s_{\Gamma})(Roo_j\T_n).\chi_{\Gamma}^s
(S(Lea_j\T_n)). \nonumber\\  \label{eq:R'Inmuts} \EEA

Expanding the formula for $S(Lea_j\T_n)$ in
terms of multiple cuts as in the previous subsection, see eq. (\ref{eq:3:S}), we get
\BEA  && [{\cal R}'{\bf\Gamma}^n(\ell(1),\ldots,\ell(n))]_{ts}=(\bar{\chi}^t_{\Gamma}-\bar{\chi}^s_{\Gamma})(\T_n)  +\sum_{l\ge 1}
(-1)^l \nonumber\\
&& \quad \sum_{j_1<\ldots<j_l} (\bar{\chi}^t_{\Gamma}-\bar{\chi}^s_{\Gamma})(Roo_{j_1}\T_n) \left\{
\prod_{m=1}^{l-1}
 \bar{\chi}^s_{\Gamma}
(Lea_{j_m}\circ Roo_{j_{m+1}}(\T_n)) \right\}  \bar{\chi}^s_{\Gamma} (Lea_{j_l}\T_n), \nonumber\\  \label{eq:2:striking} \EEA
by chopping the trunk tree $\T_n$. Finally, $\bar{\chi}_{\Gamma}^u(\T)$, $u=t$ or $s$, should be split
according to Definition \ref{def:barchi-char} (iv).

\medskip
 
Let us now make the following remark.
The difference between $[{\cal R}{\bf\Gamma}^n(\ell(1),\ldots,\ell(n))]_{ts}$ and $[{\cal R}'{\bf
\Gamma}^n(\ell(1),\ldots,\ell(n))]_{ts}$ is that
$[{\cal R}{\bf\Gamma}^n(\ell(1),\ldots,\ell(n))]_{ts}$ is obtained by first (i) splitting the measure $\mu:=\otimes_{j=1}^n
d\Gamma(\ell(j))$ into $\sum_{\sigma\in\Sigma_n}
\mu^{\sigma}\circ\sigma$ and then (ii) chopping the forests $\T_j^{\sigma}$, while $[{\cal R}'{\bf
\Gamma}^n(\ell(1),\ldots,\ell(n))]_{ts}$
is obtained by first (i) chopping the trunk tree $\T_n$ and then (ii) splitting the measures on the
trunk subtrees. Actually, as may be expected, the two operations commute.

\begin{Lemma}   \label{lem:R'=R}

$[{\cal R}'{\bf\Gamma}^n(\ell(1),\ldots,\ell(n))]_{ts}=[{\cal R}{\bf\Gamma}^n(\ell(1),\ldots,
\ell(n))]_{ts}$. Hence the regularized iterated integrals
${\cal R}{\bf\Gamma}$ satisfy the Chen and shuffle properties, and Theorem
\ref{th:2} is proved.

\end{Lemma}

{\bf Proof.}  The proof goes along the same lines as Lemma \ref{lem:barchi-char}. Let $\T_n$ be some
trunk tree with $n$ vertices and decoration $\ell$, and
$\mu:=\otimes_{j=1}^n d\Gamma(\ell(j))$.  Consider for
the moment an arbitrary character $\bar{\chi}^t_{\Gamma}$ as in Lemma \ref{lem:barchi-char}, associated
to linear forms $\phi^t_{\T}$ as in Definition \ref{def:barchi-char}. Define quite generally
\BEQ [{\cal R}_{\phi,\Gamma}(\T_n)]_{ts}:=\sum_{\sigma\in\Sigma_n} \left(
\phi^t \ast (\phi^s\circ S)\right)_{\mu^{\sigma}}(\T^{\sigma}), \EEQ
see Lemma \ref{lem:3.5}, and (see Definition \ref{def:alternative} (ii)) 
\BEQ [{\cal R'}_{\phi,\Gamma}(\T_n)]_{ts}:= \left( \chi^t_{\Gamma}\ast (\chi^s_{\Gamma}\circ S)\right)(\T_n).
\EEQ

If $\phi^t_{\T}=[{\cal R}\SkI_{\T}(\ .\ )]_t$, then ${\cal R}_{\phi,\Gamma}={\cal R}{\bf\Gamma}^n$ and
${\cal R}'_{\phi,\Gamma}={\cal R}'{\bf\Gamma}^n$. On the other hand, if $\phi^t_{\Gamma}=[\SkI_{\T}(\ .\ )]_t$,
then plainly $[{\cal R}_{\phi,\Gamma}(\T_n)]_{ts}=[{\cal R}'_{\phi,\Gamma}(\T_n)]_{ts}=
[{\bf\Gamma}^{cano,n}(\ell(1),\ldots,\ell(n))]_{ts}$, see first comment in section 2 and eq. (\ref{eq:3.12}).

Let $\sigma\in\Sigma_n$. Fix some multi-index $\vec{k}=(k_1,\ldots,k_n)$ such that $|k_1|\le \ldots\le|k_{n}|$,
and set $\mu^{\sigma}_{\vec{k}}= \tilde{\cal P}^{\{k\}}(\mu\circ\sigma)$. Then, see eq. (\ref{eq:3.7})

\BEA &&  \left( \phi^t\ast (\phi^s\circ S)\right)_{\mu^{\sigma}_{\vec{k}} }(\T^{\sigma})=
\frac{1}{|\Sigma_{\vec{k}}|} \left( \phi^t_{ \mu^{\sigma}_{\vec{k}}  } (\T^{\sigma})
+\phi^s_{ \mu^{\sigma}_{\vec{k}}  }(S(\T^{\sigma}))+ \right. \nonumber\\
&& \left.  \qquad + \sum_{\vec{v}\models V(\T^{\sigma})}
\phi^t_{Roo_{\vec{v}} ( \mu^{\sigma}_{\vec{k}} ) } (Roo_{\vec{v}}(\T^{\sigma}))
\phi^t_{Lea_{\vec{v}} ( \mu^{\sigma}_{\vec{k}} ) } (S(Lea_{\vec{v}}(\T^{\sigma}) )) \right). \EEA

Expand  $S$ according 
to eq. (\ref{eq:3:S}). This gives an expression for the $\tilde{\cal P}^{\{\vec{k}\}}$-projection of
 $[{\cal R}_{\phi,\Gamma}(\T_n)]_{ts}$.  An expression may also be obtained for the analogous regularized
quantity associated to ${\cal R}'$ by using eq. (\ref{eq:2:striking}). In the end, one gets two sums over
some subsets of $\{1,\ldots,n\}$,
\BEQ [{\cal R}_{\phi,{\cal P}^{\{\vec{k\circ\sigma}\}}\Gamma}(\T_n)]_{ts}=\sum_{J\subset \{1,\ldots,n\}} \sum_j
\phi^t_{\mu^{\sigma}_{\vec{k}}\big|_J} (\tgoth^{\sigma}_{1,J,j}) 
\phi^s_{\mu^{\sigma}_{\vec{k}}\big|_{\bar{J}} } (\tgoth^{'\sigma}_{1,\bar{J},j})  \label{eq:kR}
\EEQ
and similarly
\BEQ [{\cal R}'_{\phi,{\cal P}^{\{\vec{k\circ\sigma}\}}\Gamma}(\T_n)]_{ts}=\sum_{J\subset \{1,\ldots,n\}} \sum_j
\phi^t_{\mu^{\sigma}_{\vec{k}}\big|_J} (\tgoth^{\sigma}_{2,J,j}) 
\phi^s_{\mu^{\sigma}_{\vec{k}}\big|_{\bar{J}} } (\tgoth^{'\sigma}_{2,\bar{J},j}) \label{eq:kR'}
\EEQ

where: 
\BEQ J=V(\tgoth^{\sigma}_{1,J,j})=V(\tgoth^{\sigma}_{2,J,j}), \quad  
\bar{J}=\{1,\ldots,n\}\setminus J=
V(\tgoth^{'\sigma}_{1,\bar{J},j})
=V(\tgoth^{'\sigma}_{2,\bar{J},j}); \EEQ
\BEQ  \mu^{\sigma}_{\vec{k}}\big|_J=\otimes_{1\le j\le n, j\in J}\ 
d\tilde{\cal P}^{\{\vec{k}\circ\sigma\}} \Gamma(\ell(j)), \quad  
\mu^{\sigma}_{\vec{k}}\big|_{\bar{J}}=
\otimes_{1\le j\le n, j\in \bar{J}}\ 
d\tilde{\cal P}^{\{\vec{k}\circ\sigma\}}\Gamma(\ell(j)); \EEQ
 and $\tgoth^{\sigma}_{1,J,j}.\tgoth^{'\sigma}_{1,\bar{J},j}, 
\tgoth^{\sigma}_{2,J,j}.\tgoth^{'\sigma}_{2,\bar{J},j}\in {\cal F}^+_{n,n}$ as in the proof of Lemma
\ref{lem:barchi-char}.

In the case of the regularization scheme ${\cal R}$, each  $\tgoth_{1,J,j}^{\sigma}$ is a  forest
 such as $Roo_{\vec{v}}(\T^{\sigma}_l)$, where $\T_l^{\sigma}$ appears in the
decomposition of the permutation graph $ \T^{\sigma}$, and $\vec{v}$ is some admissible cut of $\T_l^{\sigma}$,
while $\tgoth^{'\sigma}_{1,\bar{J},j}$ is some complicated product obtained by the multiple cut decomposition
of $S(Lea_{\vec{v}}(\T_l^{\sigma}))$. In the case of ${\cal R}'$, one first splits $\T_n$ into
$(Roo_l\T_n,Lea_l\T_n)$ and then permutes the vertices of each of the two trunk subtrees, see eq.
(\ref{eq:2:striking}). 

As in the proof of Lemma \ref{lem:barchi-char}, one now proves the equality
\BEQ \sum_J \sum_j \tgoth^{\sigma}_{1,J,j}\otimes \tgoth^{'\sigma}_{1,\bar{J},j}= \sum_J
\sum_j  \tgoth^{\sigma}_{2,J,j}\otimes \tgoth^{'\sigma}_{2,\bar{J},j}  \label{eq:equality} \EEQ
by assuming that $\phi^t_{\T}=[\SkI_{\T}(\ .\ )]_t$, in which case both expressions (\ref{eq:kR}) and
(\ref{eq:kR'}) are equal.  By considering a sequence of measures $(\mu_m\circ\sigma)_{m\ge 1}$ whose
Fourier transforms converge weakly to $\del_{(\xi_1,\ldots,\xi_n)}$, one gets by Lemma \ref{lem:2:SkI}
an equation of the type
\BEA && \sum_J \left[  e^{\II (s\sum_{i\in J}\xi_i+t\sum_{i\in\bar{J}}\xi_i)}
  \sum_j F_{\tgoth^{\sigma}_{1,J,j}}( (\xi_i)_{i\in J} ) F_{\tgoth^{'\sigma}_{1,\bar{J},j}}
( (\xi_i)_{i\in\bar{J}}) \right. \nonumber\\
&&\left. \qquad  +  \sum_J    e^{\II (s\sum_{i\in J}\xi_i+t\sum_{i\in\bar{J}}\xi_i)}
\sum_j  F_{\tgoth^{\sigma}_{2,J,j}}( (\xi_i)_{i\in J} ) F_{\tgoth^{'\sigma}_{2,\bar{J},j}}
( (\xi_i)_{i\in\bar{J}})  \right] =0. \nonumber\\ \EEA
where the function $F_{\tgoth}$ has been defined in the course  of
 the proof of Lemma \ref{lem:barchi-char}.  Under the generic condition that all $\xi_J:=\sum_{i\in J}
\xi_i$, $J\subset \{1,\ldots,n\}$ are two-by-two distinct, the functions $(s,t)\mapsto
f_J(s,t):=e^{\II(s\xi_J+t\xi_{\bar{J}})}$, $J\subset\{1,\ldots,n\}$  are linearly independent. Hence, for every $J$, 
\BEQ \sum_j  F_{\tgoth^{\sigma}_{1,J,j}}( (\xi_i)_{i\in J} ) F_{\tgoth^{'\sigma}_{1,\bar{J},j}}
( (\xi_i)_{i\in\bar{J}})+F_{\tgoth^{\sigma}_{2,J,j}}( (\xi_i)_{i\in J} ) F_{\tgoth^{'\sigma}_{2,\bar{J},j}}
( (\xi_i)_{i\in\bar{J}})=0. \EEQ
By using the same arguments as in the proof of Lemma \ref{lem:barchi-char}, one obtains  eq. (\ref{eq:equality}).

\hfill \eop

%%%%%%%%%%%%%%%%%%%%%%%%%%%%%%%%%%%%%%%%%%%%%
%%%%%%%%%%%%%%%%%%%%%%%%%%%%%%%%%%%%%%%%%

\section{H\"older estimates}

%%%%%%%%%%%%%%%%%%%%%%%%%%%%%%%%%%%%%%%%%%%5
%%%%%%%%%%%%%%%%%%%%%%%%%%%%%%%%%%%%%%%%%%%%%

Let $\Gamma$ be an $\alpha$-H\"older path.
We shall now {\em choose a regularization scheme}, i.e. choose for each tree $\T$ a subset $\Z^{\T}_{reg}\subset
\Z^{\T}_+$ such that the convergence condition stated in section 2, Step 1 is verified, and
prove  that the associated regularized rough path ${\cal R}{\bf \Gamma}^n(\ell(1),\ldots,\ell(n))$
 satisfies the required H\"older properties. Following the regularization procedure as explained in section 2, one
must first  (1)
decompose  ${\cal R}{\bf \Gamma}^n(\ell(1),\ldots,\ell(n))$ into the sum over all permutations  $\sigma\in\Sigma_n$
of ${\cal R} I_{\T_j^{\sigma}} \tilde{\cal P}^+ \left( \otimes_{v\in V(\T_j^{\sigma})} \Gamma(\ell(\sigma(v)))
\right)$ as in the final step of section 2,
 and (2) show H\"older regularity with correct exponent of the {\em increment terms} 
${\cal R}\SkI_{\T}(\tilde{\cal P}^+ (\otimes_{v\in V(\T)} \Gamma(\ell(\sigma(v)))))$ and of the {\em boundary terms},
${\cal R}I_{\T}(\tilde{\cal P}^+ (\otimes_{v\in V(\T)} \Gamma(\ell(\sigma(v)))))(\partial)$, see Step 4.

%%%%%%%%%%%%%%%%%%%%%%%%%%

\subsection{Choice of the regularization scheme}

%%%%%%%%%%%%%%%%%%%%%%%%%%%%%%%%%

Recall that the whole algorithm rests on the {\em choice} of a subdomain $\Z^{\T}_{reg}\subset
\Z^{\T}_+:=\{(k_v)_{v\in V(\T)}\in \Z^{\T}\ |\ (v\twoheadrightarrow w)\Rightarrow |k_v|\ge |k_w|\}$
 for each tree $\T\in{\cal T}$. The purpose of this subsection is to propose an adequate
choice.

\begin{figure}[h]
  \centering
   \includegraphics[scale=0.35]{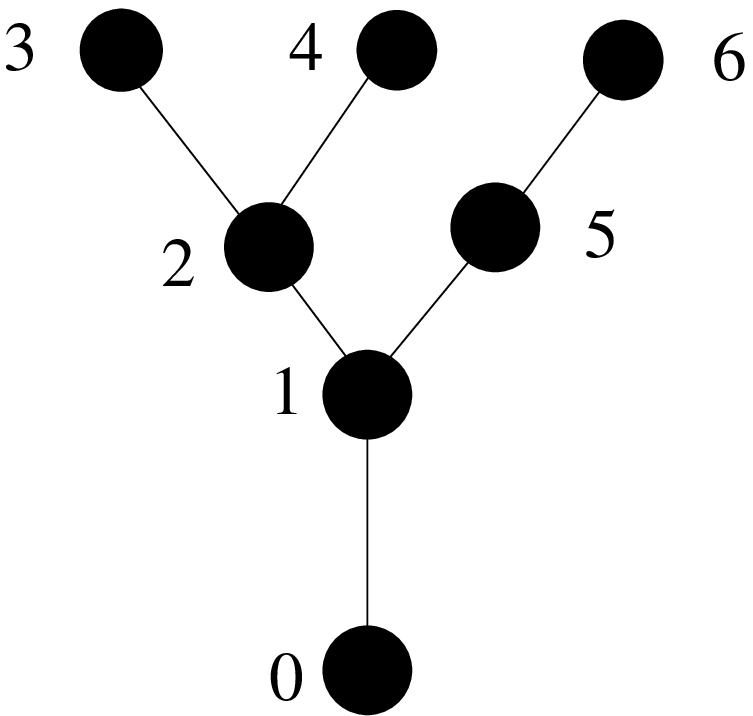}
   \caption{\small{3,4,6 are leaves; 1, 2 and 5 are nodes, 2 and 5 are uppermost; branches are e.g. 
$Br(2\twoheadrightarrow 1)=\{2\}$ or $Br(6\twoheadrightarrow 1)=\{6,5\}$; $Leaf(2)=\{3,4\}$; $w_{max}(2)=4$.}}
  \label{Fig5}
\end{figure}

We shall first need to introduce a little more terminology concerning tree structures (see Fig. \ref{Fig5}).

\begin{Definition}

Let $\T$ be a tree. 
\begin{itemize}
\item[(i)] A vertex $v$ is a {\em leaf} if no vertex connects to $v$. The set of leaves above (i.e. connecting
to) $v\in V(\T)$ is denoted by $Leaf(v)$.
\item[(ii)] Vertices at which 2 or more branches join are called {\em nodes}.
\item[(iii)] The set $Br(v_1\twoheadrightarrow v_2)$ of vertices from a leaf
or a node $v_1$ to a node $v_2$ or to the root, is called a {\em branch}
 if it does not contain any other node. By
convention, $Br(v_1\twoheadrightarrow v_2)$ includes $v_1$ and excludes $v_2$.
\item[(iv)] A node $n$ is called an {\em uppermost node} if no other node is connected to $n$.
\end{itemize}

\end{Definition}

\begin{Definition}

Let $\T$ be a tree.
If $v\in V(\T)$, we let  $w_{max}(v):=\max\{w\in V(\T)\ |\ w\twoheadrightarrow v\}$, or simply
$w_{max}(v)=v$ if $v$ is a leaf.

\end{Definition}

\begin{Definition} \label{def:2:ZTreg}

Let $\Z^{\T}_{reg}$ be the set of $V(\T)$-uples $\vec{k}=(k_v)_{v\in V(\T)}\in \Z^{\T}$ such that the following
conditions are satisfied:
\begin{itemize}
\item[(i)] if $v<w$, then $|k_v|\le |k_w|$;
\item[(ii)] if $v\in V(\T)$ and $w\in Leaf(v)$, $k_w.k_v<0$, then $|k_v|\le |k_w|-\log_2 10 - \log_2
 |V(\T)|;$
\item[(iii)] if $n\in V(\T)$ is a node, then each vertex $w\in\{w_{max}(v)\ |\ v\to n\}$
such that  $k_w.k_{w_{max}(n)}<0$ satisfies:  $|k_w|\le |k_{w_{max}(n)}|-\log_2 10-\log_2 |V(\T)|.$
\end{itemize}

\end{Definition}

\begin{Lemma} \label{lem:2:ZTreg}

Let $\vec{\xi}=(\xi_v)_{v\in\T}$ such that
 $\xi_v\in \supp(\phi_{k_v})$ for some
$\vec{k}=(k_v)_{v\in V(\T)}\in\Z^{\T}_{reg}$, where $(\phi_k)_{k\in\Z}$ is the dyadic partition of unity
defined in the Appendix.  Then, for every $v\in V$,
\BEQ |V(\T)|\ .\ |\xi_{w_{max(v)}}| \ge |\xi_v+\sum_{w\twoheadrightarrow v} \xi_w|>\half |\xi_{w_{max}(v)}|.
\EEQ

\end{Lemma}

{\bf Proof.}

The left inequality is trivial. As for the right one, 
assume first that $v$ is on a terminal branch, i.e. $Leaf(v)=\{w_{max}(v)\}$ is a singleton. Then
Definition \ref{def:2:ZTreg} (ii) implies the following: for every vertex $v'$ on the branch
between $w_{max}(v)$ and $v$, i.e. $v'\in Br(w_{max}(v)\twoheadrightarrow v)\cup \{v\}$,

-- either $\xi_{v'}$ is of the same sign as $\xi_{w_{max}(v)}$;

-- or $|\xi_{v'}|\le  \frac{|\xi_{w_{max}(v)}|}{2|V(\T)|}$
since  $|\xi_{v'}|\in (2^{|k_{v'}|-1},5\cdot 2^{|k_{v'}|-1})$ (and similarly for $|\xi_{w_{max}(v)}|$)
by the remarks following Proposition \ref{prop:1:phi}.

Hence $|\xi_v+\sum_{w\twoheadrightarrow v} \xi_w|=|\sum_{v'\in Br(w_{max}(v)\twoheadrightarrow v)\cup \{v\}}
\xi_{v'}|>\left(1-\half\frac{|\{w:w\twoheadrightarrow v\}|}{|V(\T)|} \right)
 |\xi_{w_{max}(v)}|$ and $\xi_v+\sum_{w\twoheadrightarrow v} \xi_v$ has same sign as $\xi_{w_{max}(v)}$.

\medskip

Consider now what happens at a node $n$. Let $n^+:=\{v\in V(\T)\ |\ v\to n\}$.
Assume by induction on the number of vertices that, for
all $v\in n^+$,
\BEQ \left(1+|\{w:w\twoheadrightarrow v\}|\right) |\xi_{w_{max}(v)}|\ge
|\xi_v+\sum_{w\twoheadrightarrow v} \xi_w|>\left(1-\half \frac{|\{w:w\twoheadrightarrow v\}|}{|V(\T)|} \right). \ |\xi_{w_{max}(v)}| \EEQ
and that $\xi_v+\sum_{w\twoheadrightarrow v}\xi_w$ has same sign as $\xi_{w_{max}(v)}$.
By Definition \ref{def:2:ZTreg} (iii), either $\xi_{w_{max}(v)}.\xi_{w_{max}(n)}>0$ or
$|\xi_{w_{max}(v)}|\le \frac{|\xi_{w_{max}(n)}|}{2|V(\T)|}$. Then, letting
$w_0$ be the element of $n^+$ such that $w_{max}(v_0)=w_{max}(n)$,
\BEA &&  \left(1+|\{w:w\twoheadrightarrow n\}|\right) |\xi_{w_{max}(n)}| \ge 
|\xi_n+\sum_{w\twoheadrightarrow n} \xi_w |=
 \left|\xi_n+\sum_{v\in n^+} (\xi_v+\sum_{w\twoheadrightarrow v} \xi_w) \right|  \nonumber\\
&& \qquad \qquad  \ge  \left|\xi_{v_0}+\sum_{w\twoheadrightarrow v_0} \xi_w\right|-
\left|\sum_{v\in n^+; \xi_{w_{max}}(v).\xi_{w_{max}(n)}<0} (\xi_v+\sum_{w\twoheadrightarrow v} \xi_w) \right|-|\xi_n|
\nonumber\\
&& \qquad \qquad >  \left(1-\half \frac{|\{w:w\twoheadrightarrow n\}|}{|V(\T)|} \right). \ |\xi_{w_{max}(n)}|. \EEA

\hfill \eop

%%%%%%%%%%%%%%%%%%%%%%%%%%

\subsection{A key formula for skeleton integrals}

%%%%%%%%%%%%%%%%%%%%%%%%%%

We assume in this paragraph that $\Gamma$ is smooth and denote by $\Gamma'$ its derivative. The H\"older estimates in 
subsections 4.3 and 4.4 rely on the key formula below.

\begin{Lemma} \label{lem:2:SkI}

The following formula holds:
\BEQ [\SkI_{\T}(\Gamma)]_s=(\II\sqrt{2\pi})^{-|V(\T)|} \int\ldots\int \prod_{v\in V(\T)} d\xi_v
\ .\ e^{\II s\sum_{v\in V(\T)} \xi_v} \frac{\prod_{v\in V(\T)} {\cal F}(\Gamma'(\ell(v)))(\xi_v)}{\prod_{v\in V(\T)} (\xi_v+\sum_{w\twoheadrightarrow v} \xi_w)}.\EEQ

\end{Lemma}

{\bf Proof.} We use induction on $|V(\T)|$. After stripping the root of $\T$, denoted by $0$, there remains
a forest $\T'=\T'_1\ldots\T'_J$, whose roots $0_1,\ldots,0_J$  are the vertices directly connected to $0$. Assume
\BEQ [\SkI_{\T'_j}(\Gamma)]_{x_0}=\int\ldots\int \prod_{v\in V(\T'_j)} d\xi_v \ .\ e^{\II x_0\sum_{v\in V(\T'_j)} \xi_v} F_j(\xi_{0_j},(\xi_v)_{v\in \T'_j\setminus\{0_j\}})\EEQ
for some functions $F_j$, $j=1,\ldots,J$.
Note that
\BEQ {\cal F} \left( \SkI_{\T'_j}(\Gamma) \right)(\xi_j)=\left[ \prod_{v\in V(\T'_j)\setminus\{0_j\}} \int d\xi_v 
\right]  F_j(\xi_j-\sum_{v\in V(\T'_j)\setminus\{0_j\}} \xi_v, (\xi_v)_{v\in V(\T'_j)\setminus \{0_j\} }).\EEQ

Then
\BEA
[\SkI_{\T}(\Gamma)]_s&=& \int^s d\Gamma_{x_0}(\ell(0)) \prod_{j=1}^J [\SkI_{\T'_j}(\Gamma)]_{x_0} 
\nonumber\\
&=& \frac{1}{\sqrt{2\pi}}  \int_{-\infty}^{+\infty} \frac{d\xi}{\II\xi}  e^{\II s\xi} {\cal F}\left( \Gamma'(\ell(0))
\prod_{j=1}^J \SkI_{\T'_j}(\Gamma) \right)(\xi) \nonumber\\
&=& \frac{1}{\sqrt{2\pi}} \int_{-\infty}^{+\infty} d\xi {\cal F}(\Gamma'(\ell(0)))(\xi-\sum_{j=1}^J
\xi_j) \frac{e^{\II s\xi}}{\II\xi} \ .\ \int d\xi_1\ldots\int d\xi_J   \nonumber\\
&& \qquad \qquad 
\left[ \prod_{j=1}^J \prod_{v\in V(\T'_j)\setminus \{0_j\}}   \int  d\xi_v \right]  \prod_{j=1}^J F_j(
\xi_j-\sum_{v\in V(\T'_j)\setminus\{0_j\}} \xi_v,(\xi_v)_{v\in V(\T'_j)\setminus \{0_j\} })
\nonumber\\ 
\EEA
hence the result. \hfill \eop

%%%%%%%%%%%%%%%%%%%%%%%%%%%

\subsection{Estimate for the increment term}

%%%%%%%%%%%%%%%%%%%%%%%%%%%%%%%%%%%%%

We now come back to an arbitrary $\alpha$-H\"older path and prove a H\"older estimate for the increment term, 
see eq.
(\ref{eq:incr-bdry}), which
is simply a regularized skeleton integral.

Let $\sigma\in\Sigma_n$ be a permutation, and $\T$ be one of the  forests $\T^{\sigma}_j$
appearing in the permutation graph $\T^{\sigma}$, see Lemma \ref{lem:2:sigma}. H\"older norms $|| \ . \ ||_{{\cal C}^{\gamma}}$
are defined in the Appendix. Recall $\T$ comes with a total ordering compatible with its tree partial ordering.
The $\tilde{\cal P}$-projection $\tilde{\cal P}^+$ below is defined with respect to this total ordering.

\begin{Lemma}[H\"older estimate of the  increment term]

It holds
\BEQ ||{\cal R}\SkI_{\T} \left( \tilde{\cal P}^+ (\otimes_{v\in V(\T)} \Gamma(\ell(\sigma(v)))
)\right) ||_{{\cal C}^{|V(\T)|\alpha}}<\infty.\EEQ
\label{lem:2:increment-Holder}
\end{Lemma}

{\bf Remark.} Although formal integral integrals are a priori infra-red divergent (see subsection 1.4),
the formula given in  Lemma  \ref{lem:2:SkI} for skeleton integrals delivers  infra-red convergent quantities when one
restricts the integration over 
$\xi=(\xi_v)_{v\in V(\T)}$ to the subdomain associated to $\Z^{\T}_{reg}$, see Lemma \ref{lem:2:ZTreg},
because 
\BEQ \left|\frac{{\cal F}(\Gamma'(\ell(v)))(\xi_v)}{\xi_v+\sum_{w\twoheadrightarrow v} \xi_w} \right|\lesssim
|{\cal F}(\Gamma(\ell(v)))(\xi_v)| \frac{|\xi_v|}{|\xi_{w_{max}(v)}|}\le  |{\cal F}(\Gamma(\ell(v)))(\xi_v)|
\EEQ
is bounded.

{\bf Proof.} 

We implicitly assume in the proof that $\T$ is a tree, leaving the obvious generalization to forests with
several components to the reader.

 We shall start the computations by  adapting
 the proof of a theorem in \cite{Trie}, \S 2.6.1 bounding the H\"older-Besov
norm of the product
of two H\"older functions. Write
\BEQ G(x)=\left[ {\cal R}\SkI_{\T}
 \left( \tilde{\cal P}^+ (\otimes_{v\in V(\T)} \Gamma(\ell(\sigma(v))))
\right)\right]_x.\EEQ
  By  Lemma  \ref{lem:2:SkI},
\BEA && G(x)=(\II\sqrt{2\pi})^{-|V(\T)|} \sum_{\vec{k}=(k_v)_{v\in V(\T)}\in \Z^{\T}_{reg}} 
\int_{\prod_{v\in V(\T)} \supp (\phi_{k_v})}   \prod_{v\in V(\T)} d\xi_v\ .\nonumber\\
&& \qquad \qquad .\   e^{\II x\sum_{v\in V(\T)}\xi_v}
\frac{\prod_{v\in V(\T)} {\cal F} \left( D(\phi_{k_v})\Gamma'(\ell(\sigma(v))) \right) (\xi_v)}
{\prod_{v\in V(\T)} (\xi_v+\sum_{w\twoheadrightarrow v}\xi_w)}. \nonumber\\ \label{eq:2:eq1}  
\EEA

Write, for $\vec{\xi}=(\xi_v)_{v\in V(\T)}$,
\BEQ \Theta(\vec{\xi})=\prod_{v\in V(\T)} \frac{\xi_v}
{\xi_v+\sum_{w\twoheadrightarrow v}\xi_w}\EEQ 
and 
 \BEQ \Theta_1(\vec{k})=\prod_{v\in V(\T)}
 \frac{2^{|k_v|}}{2^{|k_{w_{max}(v)}|}}.  \label{eq:Theta1} \EEQ 
Let finally  
\BEQ \Theta_{\vec{k}}(\vec{\xi}):=\prod_{v\in V(\T)}\sqrt{{\phi}_{k_v}}(\xi_v)\ .\
\frac{\Theta(\vec{\xi}) } { \Theta_1(\vec{k}) }. \EEQ
 By Lemma \ref{lem:2:ZTreg},
$||\Theta_{\vec{k}}||_{S^0(\R^{V(\T)})}$, see Proposition \ref{prop:app:Fourier-x},
  is {\em uniformly bounded} in $\vec{k}$ if $\vec{k}\in\Z^{\T}_{reg}$, which is the key
point for the following estimates.

Let $k\in\Z$. Apply the operator $D(\phi_k)$ to eq. (\ref{eq:2:eq1}): then, letting $\phi^*_{k}(\vec{\xi})
:=\phi_k(\sum_{v\in V(\T)} \xi_v)$,

\BEQ   D(\phi_k)G(x)=\left[
  \sum_{\vec{k}\in \Z^{\T}_{reg}} \Theta_1(\vec{k})D(\Theta_{\vec{k}})D(\phi^*_k) \ .\ \prod_{v\in V(\T)}
D(\sqrt{\phi_{k_v}})\Gamma(\ell(\sigma(v)))  \right](\vec{x}), \EEQ
 
where $\vec{x}=(x_v)_{v\in V(\T)}=(x,\ldots,x)$ is a vector with $|V(\T)|$ identical components.

Let $v_{max}:=\sup\{v\ |\ v\in V(\T)\}$. Note that $D(\phi^*_k)\ . \ D(\otimes_{v\in V(\T)} 
\sqrt{\phi_{k_v}})$ vanishes   except if
\BEQ \left( \sum_{v\in V(\T)} \supp (\phi_{k_v}) \right) \cap \supp (\phi_k^*)\not=\emptyset, \label{eq:2:9} \EEQ
which implies by Lemma \ref{lem:2:ZTreg}
 \BEQ |k_{v_{max}}-k|=O(\log_2 |V(\T)|); \label{eq:2:10} \EEQ
  namely, denoting by $0$ the root of $\T$,
  $|V(\T)|\ .\ |\xi_{k_{v_{max}}}|\ge |\sum_{v\in V(\T)} \xi_{k_v}|=
|\xi_{k_0}+\sum_{w\twoheadrightarrow 0}
\xi_{k_w}|> \half |\xi_{k_{v_{max}}}|$ if $\xi_v\in \supp(\phi_{k_v})$ for every $v$.

Since $\Theta_{\vec{k}},\phi^*_k \in S^0(\R^{V(\T)})$, one gets by Proposition \ref{prop:app:Fourier-x}
\BEQ ||D(\phi_k)G||_{\infty}\lesssim \sum_{\vec{k}\in\Z^{\T}_{reg},k_{v_{max}}=k}
\Theta_1(\vec{k}) \prod_{v\in V(\T)}  ||D(\sqrt{\phi_{k_v}})\Gamma(\ell(\sigma(v)))||_{\infty}.\EEQ

Since $\Gamma$ is in ${\cal C}^{\alpha}$, one obtains by Propositions \ref{prop:app:BC} and \ref{prop:app:Fourier-x}:
\BEA  ||D(\phi_k)G||_{\infty} &\lesssim& \sum_{\vec{k}\in\Z^{\T}_{reg}, k_{v_{max}}=k}
\Theta_1(\vec{k}) \prod_{v\in V(\T)} 2^{-|k_v|\alpha} \nonumber\\
&\lesssim & \sum_{\vec{k}\in\Z^{\T}_{reg},k_{v_{max}}=k}
 \prod_{v\in V(\T)} 2^{|k_v|(1-\alpha)-|k_{w_{max}(v)}|}.  \label{eq:4.19} \EEA

In other words, loosely speaking, each vertex $v\in V(\T)$ contributes a factor
$2^{|k_v|(1-\alpha)-|k_{w_{max}(v)}|}$ to $||D(\phi_k)G||_{\infty}$. If $v$ is a leaf,
then this factor is simply $2^{-|k_v|\alpha}$. Note that the upper bound
 $2^{|k_v|(1-\alpha)-|k_{w_{max}(v)}|}\le 2^{-|k_v|\alpha}$
holds true for any vertex $v$.

\medskip

Consider an uppermost node $n$, i.e. a node to which no other node is connected, together with the set
of leaves $\{w_1<\ldots<w_J\}$ above $n$, see Fig. \ref{Fig5}. Let $p_j=|V(Br(w_j\twoheadrightarrow n))|$.
On the branch number $j$,
\BEQ  2^{-|k_{w_j}|\alpha} \prod_{v\in Br(w_j\twoheadrightarrow n)\setminus\{w_j\}}  \sum_{|k_v|\le |k_{w_j}|}
  2^{|k_v|(1-\alpha)-|k_{w_j}|}  \lesssim 2^{-|k_{w_j}|\alpha p_j} \EEQ
and (summing over  $k_{w_1},\ldots,k_{w_{J-1}}$  and  over $k_{n}$)
\BEA && 2^{-|k_{w_J}|\alpha p_J}  \sum_{|k_{w_{J-1}}| \le |k_{w_J}|} 2^{-|k_{w_{J-1}}|\alpha p_{J-1}} \ \nonumber\\
&& \qquad \left( \ldots
\left( \sum_{|k_{w_1}| \le |k_{w_2}|} 2^{-|k_{w_1}|\alpha p_1}  \left( \sum_{|k_{n}| \le |k_{w_1}|}
2^{|k_n|(1-\alpha)-|k_{w_J}|} \right)\right) \ldots \right) \nonumber\\
&& \lesssim 2^{-|k_{w_J}| \alpha W(n)},\EEA
where $W(n)=p_1+\ldots+p_J+1=|\{v: v\twoheadrightarrow n\}|+1$ is the {\it weight} of $n$.

One may then consider the reduced tree $\T_n$ obtained by shrinking all vertices above $n$ (including
$n$) to {\it one} vertex  with weight $W(n)$ and perform the same operations on $\T_n$. Repeat this
inductively until $\T$ is shrunk to one point. In the end, one gets
$||D(\phi_k)G||_{\infty}
 \lesssim 2^{-|k_{v_{max}}|\alpha |V(\T)|}\lesssim 2^{-|k|\alpha|V(\T)|}$, hence $G\in {\cal C}^{|V(\T)|\alpha}.$

\hfill \eop

{\bf Remark.} Note that the above proof breaks down for the non-regularized quantitities, since
the function $\Theta_{\vec{k}}(\vec{\xi})$ is unbounded on $\Z^{\T}_+\setminus\Z^{\T}_{reg}$.
  For instance, the L\'evy area of fractional Brownian motion  diverges
below the barrier  $\alpha=1/4$, see \cite{CQ02}, \cite{Unt08}, \cite{Unt08b}.
 For deterministic, well-behaved paths $\Gamma$ with very regular,
polynomially decreasing Fourier components, the unregularized integrals are probably
well-defined at least for  $\alpha>1/2$ -- in which case the much simpler  Young integral
converges --, otherwise
the case is not even clear.

%%%%%%%%%%%%%%%%%%%%%%%%%%%%%%%%%%%%%%

\subsection{Estimate for the boundary term}

%%%%%%%%%%%%%%%%%%%%%%%%%%%%%%%%%%%%%%

We shall now prove a H\"older estimate corresponding to the boundary term. As in the previous paragraph, we let $\sigma\in
\Sigma_n$ and $\T$ be one of the  forests $\T^{\sigma}_j$, $j=1,\ldots,J_{\sigma}$. Once again, 
recall $\T$ comes with a total ordering compatible with its tree partial ordering.
The $\tilde{\cal P}$-projection $\tilde{\cal P}^+$ below is defined with respect to this total ordering.

\begin{Lemma}[H\"older regularity of the boundary term]

{\it The regularized
 boundary term $\left[{\cal R}I_{\T} \left( \tilde{\cal P}^+ (\otimes_{v\in V(\T)} \Gamma(\ell(\sigma(v))))
\right)(\partial)  \right]_{ts}$ is $|V(\T)|\alpha$-H\"older.
}
\label{lem:2:boundary-Holder}

\end{Lemma}

{\bf Proof.}

\begin{figure}[h]
  \centering
   \includegraphics[scale=0.35]{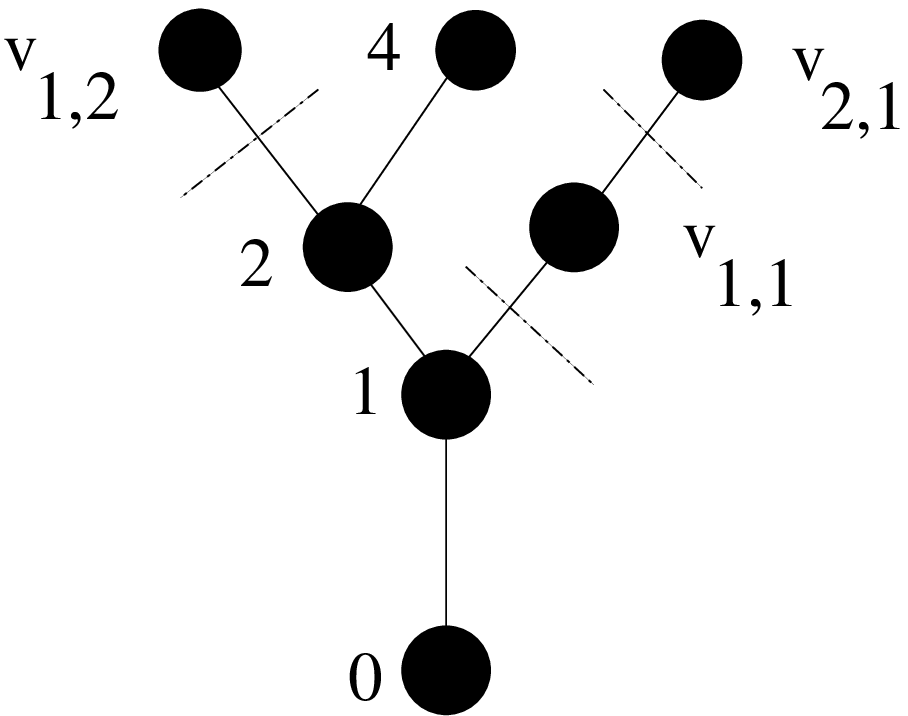}
   \caption{\small{Here $V(Roo(\T))=\{0,1,2,4\}$, $R(0)=R(4)=\emptyset$, $R(1)=\{v_{1,1}\}$,
$R(2)=\{v_{1,2}\}$.}}
  \label{Fig6}
\end{figure}

As in the previous proof, we assume implicitly that $\T$ is a tree, but the proof generalizes with only
very minor changes to the case of forests.
 Solving in terms of multiple cuts as in section 3
the recursive definition of the boundary term  
 $[{\cal R}I_{\T} \left( \tilde{\cal P}^+ (\otimes_{v\in V(\T)} \Gamma(\ell(\sigma(v)))
\right)(\partial)]_{ts}$ given in section 2, one gets  in the end  a sum of 'skeleton-type' terms of the form
(see Fig. \ref{Fig6})
\BEA &&   A_{ts}:=[\del{\cal R}\SkI_{Roo(\T)}]_{ts}
\left(  \prod_{m=1}^{l-1} [{\cal R}\SkI_{Lea_{\vec{v}_m}\circ Roo_{\vec{v}_{m+1}}
(\T)}]_s  \right) [{\cal R}\SkI_{Lea_{\vec{v}_l}(\T)}]_s
\nonumber\\
&&
\qquad \qquad  \qquad \left( \tilde{\cal P}^+ (\otimes_{v\in V(\T)} \Gamma(\ell(\sigma(v)))
\right)
 \nonumber\\ \label{eq:6:skeleton-type}
\EEA
where $\vec{v}_l=(v_{l,1}<\ldots<v_{l,J_l})\models V(\T)$, $\vec{v}_{l-1}\models V(Roo_{\vec{v}_l}\T)$,
 $\ldots$,
$\vec{v}_1=(v_{1,1},\ldots,v_{1,J_1})\models  Roo_{\vec{v}_2}(\T))$
 and one has set for short
$Roo(\T):= Roo_{\vec{v}_1}(\T)$.

\medskip

\underline{First step.}

Let $U[\vec{k}]\subset \prod_{j=1}^{J_l} \Z^{Lea_{v_{l,j}}\T}_{reg}$ such that 
 $\vec{k}=(k_{v_{l,1}},\ldots,k_{v_{l,J_l}})$ (with $|k_{v_{l,1}}|\le\ldots\le|k_{v_{l,J_l}}|$) 
 is fixed.
 Then (see after eq. (\ref{eq:4.19}) in the  proof of Lemma \ref{lem:2:increment-Holder}) each vertex $v$ contributes a factor
$2^{|k_v|(1-\alpha)-|k_{w_{max}}(v)|}\le 2^{-|k_v|\alpha}$, hence 
\BEA  || {\cal P}^{U[\vec{k}]} {\cal R}  \SkI_{ Lea_{\vec{v}_l}\T} (\otimes_{v\in V(Lea_{\vec{v}_l}\T)}
\Gamma(\ell(\sigma(v)))) ||_{\infty} &\lesssim&  
\prod_{v\in\vec{v}_l} \left[ 2^{-|k_v|\alpha}  \sum_{|k_w|\ge|k_v|,w\in Lea_v\T\setminus\{v\}} 2^{-|k_w|\alpha} \right] \nonumber\\
 & \lesssim &
\prod_{v\in\vec{v}_l}
2^{-|k_v|\alpha|V(Lea_v\T)|}.\EEA

\medskip

\underline{Second step.}

More generally, let $ B_{s}[\vec{k}]$ be the expression obtained by $\tilde{\cal P}$-projecting 
$$ \left( \prod_{m=1}^{l-1} [{\cal R}\SkI_{Lea_{\vec{v}_m}\circ Roo_{\vec{v}_{m+1}}
(\T)}]_s \right)
 [{\cal R}\SkI_{Lea_{\vec{v}_l}(\T)}]_s \left( \tilde{\cal P}^+ (\otimes_{v\in V(Lea_{\vec{v}_1}(\T))} \Gamma(\ell(\sigma(v))) )\right)$$
onto the sum of terms 
with some  fixed value of the indices  $\vec{k}=(k_{v_{1,1}},\ldots,k_{v_{1,J_1}})$.  Then
\BEQ ||B_{s}[\vec{k}] ||_{\infty} \lesssim \prod_{v\in\vec{v}_1} 2^{-|k_v|\alpha|V(Lea_{v}\T)|}
\label{eq:2:Bk} \EEQ
(proof by induction on $l$).

\medskip

\underline{Third step.}

We define
\BEA A_s(x):= && [{\cal R}\SkI_{Roo(\T)}]_x \left(
 \prod_{m=1}^{l-1} [{\cal R}\SkI_{Lea_{\vec{v}_m}\circ Roo_{\vec{v}_{m+1}}
(\T)}]_s \right)  [{\cal R}\SkI_{Lea_{\vec{v}_l}(\T)}]_s
\nonumber\\ && \qquad 
\left( \tilde{\cal P}^+ (\otimes_{v\in V(\T)} \Gamma(\ell(\sigma(v)))
\right)
\EEA
(see eq.  (\ref{eq:6:skeleton-type})), so that $A_{ts}=A_s(t)-A_s(s)$, and show that
$\sup_{s\in\R} ||x\mapsto A_s(x)||_{B^{\alpha}_{\infty,\infty}}<\infty$. Note first (see Remark following
Lemma \ref{lem:2:increment-Holder}) there is no infra-red divergence problem.

Let $V(Roo(\T))=\{w_1<\ldots<w_{max}\}$. Fix $s\in\R$ and $K\in\Z$.
By definition, and by Lemma \ref{lem:2:SkI},
\BEA  && \left( D(\phi_K)A_s \right) (x)=D(\phi_K) \left( x\mapsto
\sum_{\vec{k}=(k_{v_{1,1}},\ldots,k_{v_{1,J_1}})} \sum_{( (k_w)_{w\in V(Roo(\T))})\in S_{\vec{k}}}
\int_{\prod_{v\in V(Roo(\T))} \supp (\phi_{k_v})}  \right. \nonumber\\
&& \left.  \prod_{v\in V(Roo(\T))} d\xi_v\ .  \ e^{\II x\sum_{v\in V(Roo(\T))} \xi_v}\ 
\frac{\prod_{w\in V(Roo(\T))} {\cal F}\left( D(\phi_{k_w})\Gamma'(\ell(\sigma(w))) \right)(\xi_w)}
  {\prod_{w\in V(Roo(\T))}
(\xi_w+\sum_{w'\twoheadrightarrow w,w'\in V(Roo(\T))} \xi_{w'})}  B_{s}[\vec{k}] \right) \nonumber\\  \label{eq:4.27}
\EEA
where indices in $S_{\vec{k}}$ satisfy in particular the following conditions:

\begin{itemize}
\item[(i)]
  $|\xi_w+\sum_{w'\twoheadrightarrow w,w'\in V(Roo(\T))} \xi_{w'}|>
   \half \max\{ |\xi_{w'}|: w'\twoheadrightarrow w, w'\in V(Roo(\T))\}$ by Lemma \ref{lem:2:ZTreg};

\item[(ii)]
    $\left( \sum_{w\in V(Roo(\T))} \supp(\phi_{k_w}) \right)\cap \left( \supp(\phi_K^*)\right)\not=\emptyset$,
see eq. (\ref{eq:2:9});
\item[(iii)] for every $w\in V(Roo(\T))$, $|k_w|\le |k_{w_{max}}|$; and
\item[(iv)] for every $w\in V(Roo(\T))$, $|k_w| \le |k_v|$ for every $v\in R(w):=\{
v=v_{1,1},\ldots,v_{1,J_1}\ |\ v\to w\}$. Note that $R(w)$ may be  empty. See Fig. \ref{Fig6}.
\end{itemize}

Note that $|k_{w_{max}}-K|=O(\log_2 |V(Roo(\T))|)$ by (ii)  (see eq. (\ref{eq:2:10})).
 Hence  conditions (ii) and (iii) above are  more or
less equivalent to fixing $k_{w_{max}}\simeq K$ and letting $(k_w)_{w\in V(Roo(\T))\setminus\{w_{max}\}}$
 range
over some subset of  $[-|K|,|K|]\times\ldots\times[-|K|,|K|]$. The large fraction in eq. (\ref{eq:4.27}) contributes
to $||D(\phi_k)A_s||_{\infty}$ an overall factor bounded by $|\Theta_1(\vec{k})| \prod_{w\in V(Roo(\T))} 
2^{-|k_v|\alpha}$.

If $w\in Roo(\T)$, split $R(w)$ into $R(w)_{>}\cup R(w)_{<}$, where $R(w)_{\gtrless}
:=\{v\in R(w)\ |\ v\gtrless w_{max}\}$. Summing over indices corresponding to vertices in 
or above $R\T_{>}:=
\{v=v_{l,1},\ldots,v_{l,J_l}\ |\ v>w_{max}\}=\cup_{w\in Roo(\T)}
R(w)_{>}$, one gets by eq. (\ref{eq:2:Bk}) a quantity bounded up to a constant by
\BEQ \prod_{v\in R\T_>} \sum_{|k_v|\ge |K|} 2^{-|k_v|\alpha|V(R_v\T)|}
 \lesssim 2^{-|K|\alpha\sum_{v\in R\T_>} |V(R_v\T)|}.\EEQ

Let $w\in Roo(\T)\setminus\{w_{max}\}$ such that $R(w)_<\not=\emptyset$ 
 (note that $R(w_{max})_<=\emptyset$). Let 
 $R(w)_< = \{v_{i_1}<\ldots<v_{i_j}\}$ . Then
the sum  over $(k_v)$, $v\in R(w)_<$ contributes a factor bounded by a constant times
\BEA  && 2^{-|k_w|\alpha} \sum_{|k_{v_{i_1}}|=|k_w|}^{\infty} \sum_{|k_{v_{i_2}}|=|k_{v_{i_1}}|}^{\infty}
 \ldots \sum_{|k_{v_{i_j}}|=|k_{v_{i_{j-1}}}|}^{\infty} \nonumber\\
&& 
2^{-|k_{v_{i_1}}|\alpha |V(Lea_{v_{i_1}}\T)|} \ldots 2^{-|k_{v_{i_j}}|\alpha|V(Lea_{v_{i_j}}\T)|} 
 \lesssim 2^{-|k_w|\alpha(1+ \sum_{v\in R(w)_<} |V(Lea_v\T)|)}. \nonumber\\ \EEA

In other words, each vertex $w\in Roo(\T)$ 'behaves' as if it had a weight $1+\sum_{v\in R(w)_<} |V(R_v\T)|.$
Hence (by the same method as in the proof of Lemma \ref{lem:2:increment-Holder}), letting
$R\T_{<}:=\cup_{w\in Roo(\T)} R(w)_{<}$,
 
\BEA  ||D(\phi_K)A_s||_{\infty} &\lesssim&  2^{-|K|\alpha(|V(Roo(\T))|+\sum_{v\in R\T_<} |V(Lea_v\T)|)} \ .\ 2^{-|K|\alpha\sum_{v\in R\T_>} |V(Lea_v\T)|} \nonumber\\
&=& 2^{-|K|\alpha|V(\T)|}.\EEA

\hfill \eop

%%%%%%%%%%%%%%%%%%%%%%%%%%%%%%%%%%%%%%%%%%%%%%%%

\section{Appendix. H\"older and Besov spaces}

%%%%%%%%%%%%%%%%%%%%%%%%%%%%%%%%%%%%%%%%%%%
%%%%%%%%%%%%%%%%%%%%%%%%%%%%%%%%%%%%%%%%%%

We gather in this Appendix some definitions and technical facts about Besov spaces and H\"older norms that are required
in sections 2 and 4.

\begin{Definition}[H\"older norm]

If $f:\R^l\to\R$ is $\alpha$-H\"older continuous for some $\alpha\in(0,1)$, we let
\BEQ ||f||_{{\cal C}^{\alpha}}:=||f||_{\infty}+\sup_{x,y\in\R^l} \frac{|f(x)-f(y)|}{||x-y||^{\alpha}}.\EEQ

\end{Definition}

The space ${\cal C}^{\alpha}={\cal C}^{\alpha}(\R^l)$ of real-valued
$\alpha$-H\"older continuous functions, provided with the
above norm $|| \ ||_{{\cal C}^{\alpha}}$, is a Banach space.

\begin{Proposition} \label{prop:1:phi} \cite{Trie}

Let $l\ge 1$. There exists a family of $C^{\infty}$ functions $\phi_0,(\phi_{1,j})_{j=1,\ldots,
4^l-2^l} :\R^l\to [0,1]$,  satisfying the following conditions:

\begin{enumerate}
\item   $\supp \phi_0\subset[-2,2]^l$
and $\phi_0\big|_{[-1,1]^l}\equiv 1$.
\item Cut $[-2,2]^l$ into $4^l$ equal hypercubes of volume 1, and remove the $2^l$  hypercubes included
in $[-1,1]^l$. Let $K_1,\ldots,K_{4^l-2^l}$ be an arbitrary enumeration of the remaining 
hypercubes, and $\tilde{K}_j\supset K_j$ be the hypercube with the same center as $K_j$, but
with edges twice longer.
 Then $\supp\phi_{1,j}\subset \tilde{K}_j$, $j=1,\ldots, 4^l-2^l$.
\item Let $(\phi_{k,j})_{k\ge 2,j=1,\ldots,4^l-2^l}$ be the family of dyadic dilatations of $(\phi_{1,j})$,
namely,
\BEQ \phi_{k,j}(\xi_1,\ldots,\xi_l):=\phi_{1,j}(2^{1-k}\xi_1,\ldots,2^{1-k}\xi_l).\EEQ
Then  $(\phi_0,(\phi_{k,j})_{k\ge 1,j=1,\ldots,4^l-2^l})$ is a partition of unity subordinated to
the covering $[-2,2]^l\cup \left( \cup_{k\ge 1}\cup_{j=1}^{4^l-2^l} 2^{k-1}\tilde{K}_j\right)$, namely,
\BEQ \phi_0+\sum_{k\ge 1}\sum_{j=1}^{4^l-2^l} \phi_{k,j}\equiv 1.\EEQ
\end{enumerate}

\end{Proposition}

Constructed in this almost canonical way, the family of Fourier multipliers $(\phi_0,(\phi_{k,j}))$ is
immediately seen to be uniformly bounded for the  norm $|| \ .\ ||_{S^0(\R^l)}$ defined
in Proposition \ref{prop:app:Fourier-x} below.

If $l=1$, letting $K_1=[1,2]$ and $K_2=[-2,-1]$, we shall write $\phi_1$, resp. $\phi_{-1}$, instead
of $\phi_{1,1}$, resp. $\phi_{1,2}$, and define $\phi_k(\xi)=\phi_{\sgn(k)}(2^{1-|k|}\xi)$ for $|k|\ge 2$, so
that $\sum_{k\in\Z} \phi_k\equiv 1$ and 
\BEQ \supp \phi_0\subset[-2,2],\quad \supp \phi_k\subset [2^{k-1},5 \times 2^{k-1}],\quad
\supp \phi_{-k}\subset [-5 \times 2^{k-1},-2^{k-1}] \quad (k\ge 1). \EEQ 
In this particular case, such a family is easily constructed from an arbitrary even, 
smooth function $\phi_0:\R\to
[0,1]$ with the correct support by setting $\phi_k(\xi)={\bf 1}_{\R_+}(\xi).(\phi_0(2^{-k}\xi)-\phi_0(2^{1-k}\xi))$
and $\phi_{-k}(\xi)={\bf 1}_{\R_-}(\xi).(\phi_0(2^{-k}\xi)-\phi_0(2^{1-k}\xi))$ for every $k\ge 1$ (see
\cite{Trie2}, \S 1.3.3).

\medskip

In order to avoid setting apart the one-dimensional case, we let $\I_l:=\Z$ if $l=1$, and
$\I_l=\{0\}\cup\{(k,j)\ |\ k\ge 1, 1\le j\le 4^l-2^l\}$ if $l\ge 2$. Also, if $l\ge 2$, we
define $|\kappa|=k\ge 1$ if $\kappa=(k,j)$ with $k\ge 1$.

\begin{Definition} \label{def:phitildek}
Let $(\tilde{\phi}_{\kappa})_{\kappa\in \I_l}$ be the partition of unity of $\R^l$, $l\ge 1$ defined by
(see Proposition \ref{prop:1:phi}) : 
\begin{itemize}
\item[(i)]
\BEQ \tilde{\phi}_0:={\bf 1}_{[-1,1]^l},\quad \tilde{\phi}_{1,j}:={\bf 1}_{K_j};\EEQ
\item[(ii)] if $k\ge 2$,
\BEQ \tilde{\phi}_{k,j}(\xi_1,\ldots,\xi_l):=\tilde{\phi}_{1,j}(2^{1-k}\xi_1,\ldots,2^{1-k}\xi_l).\EEQ
\end{itemize}
\end{Definition}

We use this auxiliary partition several times in the text.

\begin{Definition} \cite{Trie}

Let $\ell_{\infty}(L_{\infty})$ be the space of sequences $(f_{\kappa})_{\kappa\in\I_l}$ of a.s. bounded
functions $f_{\kappa}\in L_{\infty}(\R^l)$ such that 
\BEQ ||f_{\kappa}||_{\ell_{\infty}(L_{\infty})}:=\sup_{\kappa\in\I_l} ||f_{\kappa}||_{\infty}<\infty.\EEQ
\end{Definition}

Let ${\cal S}'(\R^l,\R)$ be the dual of the Schwartz space of rapidly decreasing functions on $\R^l$. As
well-known, it includes the space of infinitely differentiable slowly growing functions.

The following definition is classical. Recall that the Fourier transform $\cal F$ has been defined at the end of
the Introduction.

\begin{Definition}[Fourier multipliers]

Let $m:\R^l\to\R$ be an infinitely differentiable slowly growing function. Then 
\BEQ D(m): {\cal S}'(\R^l,\R)\to{\cal S}'(\R^l,\R),\quad \phi\mapsto {\cal F}^{-1}(m\cdot {\cal F}\phi) \EEQ
defines a continuous operator.

\end{Definition}

In other words, $m$ is a Fourier multiplier of ${\cal S}'(\R^l,\R)$.

\begin{Definition} \cite{Trie}

Let $B^{\alpha}_{\infty,\infty}(\R^l):=\{f\in {\cal S}'(\R^l,\R)\ |\ ||f||_{B^{\alpha}_{\infty,\infty}}<\infty\}$
where
\BEA
 ||f||_{B^{\alpha}_{\infty,\infty}} &:=& ||2^{\alpha|\kappa|}D(\phi_{\kappa})f||_{\ell_{\infty}(L_{\infty})}
 \nonumber\\
&=& \sup_{\kappa\in\I_l} 2^{\alpha|\kappa|} ||D(\phi_{\kappa})f||_{\infty}.\EEA

\end{Definition}

\begin{Proposition} (see \cite{Trie}, \S 2.2.9) \label{prop:app:BC}

For every $\alpha\in(0,1)$, $B^{\alpha}_{\infty,\infty}(\R^l)={\cal C}^{\alpha}(\R^l)$, and the two norms
$||\ ||_{{\cal C}^{\alpha}}$ and $||\ ||_{B^{\alpha}_{\infty,\infty}}$ are equivalent.

\end{Proposition}

We shall sometimes call $||\ ||_{B^{\alpha}_{\infty,\infty}}$ the {\it H\"older-Besov norm}.

 Let us finally give a criterion for a function $m$ to be a Fourier multipliers
of the Besov space $B_{\infty,\infty}^{\alpha}$:

%\begin{Proposition}[maximal inequality] (see \cite{Trie}, \S 2.1.2 p. 29)

%Let, for  $a>0$ and $f\in{\cal S}'$,
%\BEQ (\phi^*_k f)(x):=\sup_{y\in\R} \frac{|  (D(\phi_k)f)(x-y)|}{1+(2^k|y|)^a}.\EEQ

%Then 
%\BEQ || 2^{\alpha k} \phi^*_k f ||_{\ell_{\infty}(L_{\infty})}\le C ||f||_{B^{\alpha}_{\infty,\infty}}.\EEQ
%\end{Proposition}

\begin{Proposition}[Fourier multipliers](see \cite{Trie}, \S 2.1.3, p.30)

Let $\alpha\in(0,1)$ and $m:\R^l\to\R$ be an infinitely differentiable  function such that
\BEQ ||m||_{S^0(\R^l)}:=\sup_{|j|\le l+5} \sup_{\vec{\xi}\in\R^l} |(1+||\vec{\xi}||)^{|j|}  m^{(j)}(\vec{\xi})|
 <\infty\EEQ
where $j=(j_1,\ldots,j_l)$, $|j|=j_1+\ldots+j_l$ and $m^{(j)}:=\partial_{\xi_1}^{j_1}\ldots\partial_{\xi_l}^{j_l}m.$
 Then there exists a constant $C$ depending only on $\alpha$, such that
\BEQ ||D(m)f||_{B^{\alpha}_{\infty,\infty}} \le C ||m||_{S^0(\R^l)} ||f||_{B^{\alpha}_{\infty,\infty}}.\EEQ

\label{prop:app:Fourier-x}
\end{Proposition}

The space $S^0(\R^l)$ contains the space of translation-invariant pseudo-differential symbols
of order $0$ (see for instance \cite{BJR}, Definition 1.1, or \cite{Treves}).

%\begin{Definition}

%The multiplicative property implies in particular the following identity for the
% twice iterated integral ${\cal A}_{ts}:=\int_s^t d\Gamma_{x_1}(1)\int_s^{x_1}d\Gamma_{x_2}(2)$ 
%(which is a way to  measure the area generated by the first two components of $\Gamma$) :
%\BEQ {\cal A}_{ts}={\cal A}_{tu}+{\cal A}_{us}+ (B_t(1)-B_u(1))(B_u(2)-B_s(2)) \label{eq:0:mult} \EEQ
%while the geometric property implies
%\BEA && \int_s^t dB_{t_1}(1)\int_s^{t_1} dB_{t_2}(2)+\int_s^t dB_{t_2}(2)\int_s^{t_2} dB_{t_1}(1) \nonumber\\
%&&=
%\left(\int_s^t dB_{t_1}(1)\right) \left(\int_s^t dB_{t_2}(2)\right)=(B_t(1)-B_s(1))(B_t(2)-B_s(2)). \nonumber\\
%\EEA

%%%%%%%%%%%%%%%%%%%%%%%%%%%%%%%%%%%%%%%%%%%
%%%%%%%%%%%%%%%%%%%%%%%%%%%%%%%%%%%%%%%%%%%

\end{document}